\input amstex
\input epsf
\documentstyle{amsppt}
\magnification=1200
\hcorrection{0.2in}
\vcorrection{-0.1in}
\define\xuparrow#1{%
{\left\uparrow\vbox to #1{}\right.\kern-\nulldelimiterspace}
}
\NoRunningHeads
\NoBlackBoxes
\topmatter
\title Finite Quotient of Join in Alexandrov Geometry
\endtitle
\author Xiaochun Rong \footnote{Supported partially by
NSF Grant DMS 1106517 and a research found from Capital normal university.
\hfill{$\,$}} \& Yusheng Wang \footnote{Supported partially by
NFSC 11471039.\hfill{$\,$}}
\endauthor
\abstract  Given two $n_i$-dimensional Alexandrov spaces $X_i$ of curvature
$\ge 1$, the join of $X_1$ and $X_2$ is an $(n_1+n_2+1)$-dimensional Alexandrov
space $X$ of curvature $\ge 1$, which contains $X_i$ as convex subsets
such that their points are $\frac \pi2$ apart.  If a group acts isometrically
on a join  that preserves $X_i$, then the orbit space is called quotient
of join. We show that an $n$-dimensional Alexandrov space $X$
with curvature $\ge 1$ is isometric to a finite quotient of join, if
$X$ contains two compact convex subsets $X_i$ without boundary such
that $X_1$ and $X_2$ are at least $\frac \pi2$ apart and
$\dim(X_1)+\dim(X_2)=n-1$.
\endabstract
\address Mathematics Department, Capital Normal University,
Beijing, P.R.C. \newline.\hskip3mm
Mathematics Department, Rutgers University, New Brunswick,
NJ 08903, U.S.A
\endaddress
\email rong\@math.rutgers.edu
\endemail
\address School of Mathematical Sciences (and Lab. math. Com.
Sys.), Beijing Normal University, Beijing, 100875 P.R.C.
\endaddress
\email wyusheng\@bnu.edu.cn
\endemail

\endtopmatter
\document

\vskip2mm

\head 0. Introduction
\endhead

\vskip4mm

Consider a compact Riemannian manifold $M$ of positive sectional
curvature, normalized to $\text{sec}_M\ge 1$. The diameter,
$\text{diam}(M)\le \pi$ (Bonnet theorem), and ``$=$'' holds if and
only if $M$ is isometric to the standard sphere ([Ch]). If
$\text{diam}(M)>\frac \pi2$, then $M$ is homeomorphic to the
standard sphere ([GS]). If $\text{diam}(M)=\frac \pi2$, then the
Diameter rigidity theorem ([GG]) says that $M$ is either
homeomorphic to a sphere  or locally isometric to a rank one
symmetric space with the standard metric (cf. [Wi]).

In this paper, we will explore a rigidity of a finite quotient of join in
Alexandrov geometry (see Theorem A); which is a necessary step toward a
classification for Alexandrov spaces of curvature $\ge 1$ and
diameter $\frac \pi2$. Furthermore, understanding such a finite quotient
of join is crucial to solve the Soul Conjecture in Alexandrov geometry
by Perel'man ([Pe], [Li], [RW]).

Alexandrov geometry was introduced by Burago-Gromov-Perelman in
[BGP], and extensive study has been done since then.
%(cf. [Ka], [Pet]).
An Alexandrov space with curvature $\ge \kappa$  is a locally
compact length metric space such that every geodesic triangle looks
fatter than a corresponding triangle in the simply connected
$2$-space form of constant curvature $\kappa$, i.e., Toponogov's
comparison theorem holds (see Theorem 1.1 below).
%A partial motivation for
%studying Alexandrov spaces is that Gromov-Hausdorff limits of
%Riemannian $n$-manifolds with sectional curvature uniformly bounded
%by $\kappa$ may not be a Riemannian manifold but an Alexandrov space
%with curvature $\ge \kappa$.

Let $\text{Alex}^n(\kappa)$ denote the set of all complete
$n$-dimensional Alexandrov spaces with $\text{cur}\ge \kappa$. Many
results in Riemannian geometry based on Toponogov's comparison
theorem has been generalized to Alexandrov spaces ([BGP], cf.
[AKP]). Concerning a `classification' of $X\in \text{Alex}^n(1)$
with $\text{diam}(X)\ge \frac \pi2$, the following are known: $X\in
\text{Alex}^n(1)$ satisfies $\text{diam}(X)\le \pi$, ``$=$'' holds
if only only if $X$ is isometric to spherical suspension,
$S(\Sigma)$ with $\Sigma\in \text{Alex}^{n-1}(1)$, and if
$\text{diam}(X)>\frac \pi2$, then $X$ is homeomorphic to some
$S(\Sigma)$ ([BGP]). Note that here $\Sigma\in \text{Alex}^{n-1}(1)$
($n\ge 4$) can have infinitely many possible topological types, and
thus unlike a complete classification in Riemannian case, the best
`classification' one can expect is to find a rigid underlying
geometric structure involving Alexandrov spaces of lower dimension.

Guided by the Riemannian situation, a challenging problem is to
explore a rigid structure on $X\in\text{Alex}^n(1)$ with
$\text{diam}(X)=\frac \pi2$. Similar to the Riemannian case ([GG]),
$X$ always has two compact convex subsets which are $\frac
\pi2$-apart: let $p\in X$, then $X_1=\{x\in X |\, |xp|=\frac \pi2\}$
and $X_2=\{x\in X|\, |xX_1|=\frac \pi2\}$ are compact convex subsets
\footnote{We say that a subset $A\subseteq X$ is convex (resp.
totally convex) if for any $x, y\in A$, there is a minimal geodesic
(resp. all minimal geodesics) jointing $x$ with $y$ that is in
$A$.}, where `$|\cdot \cdot|$' denotes the distance. Let
$m=\dim(X_1)+\dim(X_2)$. Roughly, our main result asserts that $m\le
n-1$ and ``$=$'' implies a rigid underlying geometric structure, and
thus the remaining `classification' reduces to $m\le n-2$ (see
Remark 0.5).

Let's now turn to a known construction for Alexandrov $n$-spaces of
curvature $\ge 1$ ([BGP]). Let $Y_i\in \text{Alex}^{n_i}(1)$ ($i=1,
2$), and $Y_i=\{p_i\}$ or $\{p_i,q_i\}$ with $|p_iq_i|=\pi$ if
$n_i=0$. Then the following space
$$Y_1*Y_2=Y_1\times Y_2\times \left[0,\frac \pi2\right]/(y_1,y_2,0)\sim (y_1,y_2',0),
\left(y_1,y_2,\frac \pi2\right)\sim \left(y_1',y_2,\frac
\pi2\right),$$ equipped with the metric
$$\cos |(y_1,y_2,t)(y_1',y_2',s)|=\cos s\cos t\cos |y_1y_1'|+\sin s\sin t
\cos |y_2y_2'|,$$ is called the join of $Y_1$ and $Y_2$. Especially,
if $Y_1=\{p,q\}$ with $|pq|=\pi$, then $Y_1*Y_2$ is also called a
spherical suspension over $Y_2$, $S(Y_2)$. If $Y_1=\{p\}$, then the
join is a half spherical suspension over $Y_2$, $S_+(Y_2)$.

It is easy to see that $Y_1*Y_2\in \text{Alex}^n(1)$ with the
following properties:

\noindent (0.1.1) $|p_1p_2|=\frac \pi2$ for all $p_i\in Y_i$.

\noindent (0.1.2) $n_1+n_2=n-1$.

\noindent (0.1.3) For $p_i\in Y_i$, there is a unique minimal
geodesic jointing $p_1$ with $p_2$.

Observe that $Y_1*Y_2$ is a
Riemannian manifold if and only if $Y_i$ and $Y_1*Y_2$ are unit spheres.

Assume that a group $\Gamma$ acts effectively on $Y_1$ and $Y_2$ by isometries.
By (0.1.3), the $\Gamma$-action uniquely extends to an isometric
$\Gamma$-action on $Y_1*Y_2$. The quotient space,
$(Y_1*Y_2)/\Gamma\in \text{Alex}^n(1)$ satisfies (0.1.1) and
(0.1.2), but not (0.1.3). We will call $(Y_1*Y_2)/\Gamma$ a quotient
of join. From the construction of a quotient of join, we observe the
following properties:

\noindent (0.1.4) $\text{diam}((Y_1*Y_2)/\Gamma)=\max \{\frac \pi2,\text{diam}(Y_i/\Gamma)\}$,
and $\dim((Y_1*Y_2)/\Gamma)=\frac \pi2$ if and only if $\text{diam}(Y_i/\Gamma)\le \frac \pi2$.

\noindent (0.1.5) If both $Y_1$ and $Y_2$ have an empty boundary \footnote{If $\dim(Y_i)=0$, then that
$Y_i$ has an empty boundary means that $Y_i=\{p_i,q_i\}$ with
$|p_iq_i|=\pi$.}, then $(Y_1*Y_2)/\Gamma$
has an empty boundary. Note that if $\Gamma$ acts trivially
on $Y_1$, then the joint structure on $Y_1*Y_2$ descends to a join structure on
$Y_1*(Y_2/\Gamma)$, which may have an non-empty boundary.

The main result in this paper is the following rigidity of a finite quotient of
join:

\vskip2mm

\proclaim{Theorem A} Let $X\in \text{Alex}^n(1)$, and let $X_1, X_2$
be two compact convex subsets in $X$ such that
$|X_1X_2|\triangleq\min\{|x_1x_2|\ |\, x_i\in X_i\}\geq \frac \pi2$.
Then

\noindent (A1) $\dim(X_1)+\dim(X_2)\le n-1$, provided either $X_1$
or $X_2$ has an empty boundary.

\noindent (A2) If both $X_1$ and $X_2$ have an empty boundary and
$\dim(X_1)+\dim(X_2)=n-1$, then $X$ is isometric to a finite
quotient of join (which is a join when $\Gamma=e$). Precisely, there is $p_i\in X_i$ ($i=1, 2$) and a
finite group $\Gamma$ acting effectively and isometrically on
$(\Sigma_{p_i}X_i)^\perp$ such that
$$X\overset{\text{isom}}\to\cong [(\Sigma_{p_1}X_1)^\perp * (\Sigma_{p_2}X_2)^\perp]
/\Gamma,\hskip2mm X_1\overset{\text{isom}}\to
\cong(\Sigma_{p_2}X_2)^\perp/\Gamma,\hskip2mm X_2
\overset{\text{isom}}\to \cong(\Sigma_{p_1}X_1)^\perp/\Gamma,$$
where $\Sigma_qX$ denotes the space of directions of $X$ at $q$, and
$(\Sigma_{p_i}X_i)^\perp=\{u\in \Sigma_{p_i}X |\, |uv|=\frac \pi2,\,
v\in \Sigma_{p_i}X_i\}$.
\endproclaim

Note that in Theorem A there is no requirement of $\text{diam}(X)=
\frac \pi2$ (see (0.1.4)), and $X$ satisfying (A2) has an empty boundary
(see (0.1.5)). One may review (A1) as an analog in Alexandrov geometry of
Proposition 1.4 in [GG], where $X$ is a Riemannian manifold with
$\sec_X\geq 1$ and $\text{diam}(X)=\frac\pi2$.

In Riemannian case, i.e., $X$ is a Riemannian manifold, Theorem A
has the following corollary.

\proclaim{Corollary 0.2} Let $M$ be a compact $n$-manifold with
$\text{sec}_M\ge 1$, and let $M_i\subset M$ be two compact totally
geodesic submanifolds with $|M_1M_2|\ge \frac \pi2$. Then
$\dim(M_1)+\dim(M_2)\leq n-1$, and ``='' implies that $M$ is
isometric to a spherical space form with a quotient join structure.
\endproclaim

Here, if $\dim(M_i)=0$, then $M_i$ is a single point. In Corollary 0.2,
``$\dim(M_1)+\dim(M_2)\leq n-1$'' is also seen by
Frankel's theorem ([Fr]). Indeed, Corollary 0.2 recovers the spherical
space form case in the Diameter rigidity theorem in [GG]. Moreover,
restricted to Riemmanian case, our approach is different from [GG];
ours relies on finding the underlying join structure (Remark 5.1).

Note that not every spherical space form satisfies Corollary 0.2 (e.g.,
any spherical $3$-space form whose  fundamental group is not
cyclic). However, based on [Ro] we show that up to a finite normal
covering space of a uniform bounded order, all spherical space forms
satisfy Corollary 0.2.

\proclaim{Corollary 0.3} There is a constant $w(n)>0$ such that any
spherical $n$-space form has a normal covering space of order $\le
w(n)$ which satisfies Corollary 0.2.
\endproclaim

Let's make a few remarks on the above results.

\remark{Remark \rm 0.4} In the proof (A2), we construct a finite
group $\Gamma$ acting isometrically on $(\Sigma_{p_i}X_i)^\perp$ at
some $p_i\in X_i$, which extends to a unique isometric action on the
join of $(\Sigma_{p_1}X_1)^\perp*(\Sigma_{p_2}X_2)^\perp$ with orbit
space isometric to $X$. Recall that given an orbit space $Y/\Gamma$
(with a stratification and assigned isotropy groups), it may not be
possible to recover the $\Gamma$-action on $Y$. Here that we are
able to recover the $\Gamma$-action is based on (0.1.1) (see Lemma
1.2 below) and the rigidity part of Toponogov's comparison theorem
which implies that everywhere there is isometrically embedded
spherical triangle; precisely, for all $p_i, p_i' \in X_2$, the
hinge, $([p_1p_2],[p_1p_1'],\frac \pi2)$ or $([p_2p_1],[p_2p_2'],
\frac \pi2)$ bounds an isometrically embedded spherical triangle
([GM]), where $[pq]$ denotes a minimal geodesic from $p$ to $q$.

We notice that on compact manifolds of non-positive sectional
curvature, the higher rank rigidity theorem (cf. [BGS], and reference
within) is based on the geometry that at every point and in every
direction, there is isometrically embedded plane.
\endremark

\vskip2mm

\remark{Remark \rm 0.5} We mention that recently, an isometric
classification has been obtained by [SSW] for Riemannian
$n$-manifolds of $\text{sec}\ge 1$ which contains two compact
totally geodesic submanifolds $X_i$ such that $|X_1X_2|\ge \frac
\pi2$ and $\dim(X_1)+\dim(X_2)=n-2$.
\endremark

\vskip2mm

Given a subset $A\subset X\in
\text{Alex}^n(1)$, let $A^{\ge\frac \pi2}=\{x\in X|\ |xA|\geq \frac\pi2\}$
and let $A^{=\frac\pi2}=\{x\in X|\ |xa|=\frac\pi2, \forall\ a\in A\}$.
Note that $A^{\ge\frac \pi2}$ is totally convex in $X$ (see (1.1.1)).

Our last theorem asserts that the conclusion of (A2) still holds if one of
$X_i$ has a non-empty boundary and some addition restrictions.

\proclaim{Theorem B} Let $X\in \text{Alex}^n(1)$,
and let $X_1, X_2\subset X$ be two compact convex subsets. Assume
that $|X_1X_2|\geq \frac \pi2$, $\partial X_1=\emptyset$, $\partial
X_2\neq\emptyset$, and $\dim(X_1)+\dim(X_2)=n-1$. If $X_2^{\geq\frac
\pi2}=X_2^{=\frac \pi2}$ and $(\Sigma_{p_2}X_2)^{\ge\frac \pi2}
=(\Sigma_{p_2}X_2)^{=\frac \pi2}$ for all $p_2\in \partial X_2$,
then $X$ is isometric to a finite quotient of join as in Theorem A.
\endproclaim

In Theorem B, if $X_2=\{p_2\}$, then $X$ is either isometric to $S_+(X_1)$,
or  isometric to $S(\Sigma_{p_2}X)/\Bbb Z_2$
and $X_1$ is isometric to $(\Sigma_{p_2}X)/\Bbb Z_2$.

\vskip2mm

We now give a brief description on our approach to Theorem A.

We first give some conventions on the notations in the paper. Let
$\uparrow_{x_1}^{x_2}$ denote the direction at $x_1$ of a given
$[x_1x_2]$, let $\Uparrow_{x_1}^{x_2}=\{(\uparrow_{x_1}^{x_2})_i\}$
denote the collection of $\uparrow_{x_1}^{x_2}$'s, and
let $\Uparrow_{x_1}^{Z} =\bigcup_{z\in Z}\Uparrow_{x_1}^z$.

Observe that $|X_1X_2|\ge \frac \pi2$ implies that for all $p_i\in
X_i$, $|p_1p_2|=\frac \pi2$ (see Lemma 1.3 below), and thus
$\Uparrow_{p_1}^{X_2}\subseteq (\Sigma_{p_1}X_1)^\perp$. We consider
the multi-value map $f_{p_1}: X_2\to (\Sigma_{p_1}X_1)^\perp$
defined by $f_{p_1}(p_2')=\Uparrow_{p_1}^{p_2'}$. Then,
$f_{p_1}(X_2)=\Uparrow_{p_1}^{X_2}$. Based on the rigidity part of
Toponogov's comparison theorem, a crucial observation is:

\vskip2mm

\noindent (0.6) $f_{p_1}$ is a `radial isometry', i.e., any given
$[p_1p_2]$ and $[p_2p_2']$ bounds an embedded standard spherical
triangle, whose third side, $[p_1p_2']$, determines
$[\uparrow_{p_1}^{p_2}\uparrow_{p_1}^{p_2'}]\subset
(\Sigma_{p_1}X_1)^\perp$ such that $f_{p_1}: [p_2p_2']\to
[\uparrow_{p_1}^{p_2}\uparrow_{p_1}^{p_2'}]$ is an isometry.

\vskip2mm

By Toponogov's comparison theorem and (0.6), $f_{p_1}$ is distance
non-decreasing, so by induction $\dim(X_1)+\dim(X_2)\le
\dim(X_1)+\dim((\Sigma_{p_1}X_1)^\perp)\leq n-1$ ((A1)).

In the proof of (A2), a key is to show that
$\Uparrow_{p_1}^{X_2}=(\Sigma_{p_1}X_1)^\perp$ and (and
$\Uparrow_{p_2}^{X_1}=(\Sigma_{p_2}X_2)^\perp$) (see Key Lemma 1.7),
and it suffices to show that $\Uparrow_{p_1}^{X_2}$ ($\subseteq
(\Sigma_{p_1}X_1)^\perp$) is a compact top dimensional locally
convex subset without boundary. Using (0.6), at $p_2$ and any
$\uparrow_{p_1}^{p_2}\in f_{p_1}(p_2)$, we define a natural
multi-valued `tangent map' of $f_{p_1}$,
$\text{D}f_{p_1}:\Sigma_{p_2}X_2\to
\Sigma_{\uparrow_{p_1}^{p_2}}(\Sigma_{p_1}X_1)^\perp$, that is also
distance non-decreasing. More importantly, $\text{D}f_{p_1}$ is also
a radial isometry as well, which allows us to proceed an induction
to prove that $\Uparrow_{p_1}^{X_2}=(\Sigma_{p_1}X_1)^\perp$.

In practice, we first show that for all $p_i\in X_i$,
$|\Uparrow_{p_1}^{p_2}|\le m<\infty$ and pairs of points at which
$|\Uparrow_{p_1}^{p_2}|=m$ are dense in $X_i$ (see Key Lemma 1.6).
The boundedness is used to show that $\text{D}f_{p_1}$ is a radial
isometry and $|\text{D}f_{p_1}(\xi)|$ $(\xi\in \Sigma_{p_2}X_2)$ has
a similar boundedness.

Once $\Uparrow_{p_1}^{X_2}=(\Sigma_{p_1}X_1)^\perp$ is established,
fixing $p_i$ with $|\Uparrow_{p_1}^{p_2}|=m$ we can construct a
finite group $\Gamma_2$ acting isometrically on
$(\Sigma_{p_1}X_1)^\perp$ such that
$\Gamma_2(\uparrow_{p_1}^{p_2})=\Uparrow_{p_1}^{p_2}$,  and
$(\Sigma_{p_1}X_1)^\perp/\Gamma_2$ is isometric to $X_2$. Given
$\uparrow_{p_2}^{p_1}$ and $\uparrow_{p_2}^{p_1}\ne
(\uparrow_{p_2}^{p_1})_i\in \Uparrow_{p_2}^{p_1}$, a geodesic
jointing the two points determines a geodesic loop $\gamma_i\subset
X_1$ at $p_1$. For each $\gamma_i$, we define an isometry on
$(\Sigma_{p_1}X_1)^\perp$ via a successive `parallel transport'
along $\gamma_i$ piece by piece. The independence of the isometry on
the choice of $\gamma_i$ relies on the fact that isometrically
embedded spherical triangles exist everywhere (see Remark 0.4).

Similarly, we can construct a finite group $\Gamma_1$ acting
isometrically on $(\Sigma_{p_2}X_2)^\perp$. And from the
construction of $\Gamma_i$, it is not hard to show that $\Gamma_1$
is isomorphic to $\Gamma_2$.  Consequently, we obtain an isometric
$\Gamma$-action on
$(\Sigma_{p_1}X_1)^\perp*(\Sigma_{p_2}X_2)^\perp$, and it remains to
show that $X$ is isometric to the orbit space. Note that a
consequence of $\Uparrow_{p_1}^{X_2}=(\Sigma_{p_1}X_1)^\perp$ is
that for any $x\in X-(X_1\cup X_2)$, there are unique $p_i\in X_i$
and a minimal geodesic jointing $p_1$ with $p_2$ such that $x\in
[p_1p_2]$. This property guarantees that the natural map from
$[(\Sigma_{p_1}X_1)^\perp*(\Sigma_{p_2}X_2)^\perp]/\Gamma$ to $X$ is
a bijection, which preserves the distance (this also relies on the
rigidity part of Toponogov's comparison theorem).

\vskip4mm

\head 1. The Geometry of $\frac \pi2$-apart Convex Subsets
\endhead

\vskip4mm

In the proof of Theorem A, we will often apply Toponogov's comparison theorem,
in particular its rigidity part in (1.1.3). Recall that $|pq|$ (resp. $[pq]$)
denotes the distance (resp. a minimal geodesic) between $p$ and $q$.

\vskip2mm

\proclaim{Theorem 1.1 ([BGP])} Let $X\in\text{Alex}^n(\kappa)$, and
let $\Bbb S^2_\kappa$ be the complete, simply connected $2$-manifold
of curvature $\kappa$.

\noindent {\rm(1.1.1)} To any $p\in X$ and $[qr]\subset X$,
we associate $\tilde p$ and a $[\tilde q\tilde r]$
in $\Bbb S^2_\kappa$ with $|\tilde p\tilde q|=|pq|,|\tilde p\tilde
r|=|pr|$ and $|\tilde r\tilde q|=|rq|$. Then for any $s\in[qr]$ and
$\tilde s\in[\tilde q\tilde r]$ with $|qs|=|\tilde q\tilde s|$,
we have that $|ps|\geq|\tilde p\tilde s|$.

\vskip1mm

\noindent {\rm(1.1.2)} To any $[qp]$ and $[qr]$ in $X$,
we associate $[\tilde q\tilde p]$ and $[\tilde q\tilde r]$
in $\Bbb S^2_\kappa$ with $|\tilde q\tilde p|=|qp|$, $|\tilde q\tilde r|=|qr|$
and $\angle\tilde p\tilde q\tilde r=\angle pqr$. Then we have that
$|\tilde p\tilde r|\geq|pr|$.

\vskip1mm

\noindent {\rm(1.1.3) \bf([GM])} If equality in {\rm (1.1.2)}
(resp. in {\rm (1.1.1)} for some interior point $s$ in $[qr]$)
holds, then there exists a $[pr]$ (resp. $[qp]$ and $[pr]$) such
that $[qp]$, $[qr]$ and $[pr]$ bounds a surface which is convex and
can be isometrically embedded into $\Bbb S^2_\kappa$.
\endproclaim

We first observe the following fact on $X_1$ and $X_2$ in Theorem A.

\proclaim{Lemma 1.2} Let $X\in \text{Alex}^n(1)$. Assume that $X_1,
X_2$ are compact convex subsets such that $|X_1X_2|\ge \frac \pi2$.
If $\partial X_1=\emptyset$ (or $\partial X_2=\emptyset$), then
$|x_1x_2|=\frac\pi2$ for all $x_i\in X_i$.
\endproclaim

Lemma 1.2 is a consequence of the following lemma.

\proclaim{Lemma 1.3 ([Ya])} Let $X\in {Alex}^n(1)$, and let $A$ be a
compact (locally) convex subset in $X$. If $\partial A=\emptyset$,
then $A^{\geq \frac{\pi}{2}}=A^{=\frac{\pi}{2}}$.
\endproclaim

Note that using (1.1.2), one can derive Lemma
1.3 by induction on $n$.

From now on, we let $X_1, X_2\subset X\in \text{Alex}^n(1)$ be
compact convex subsets with $$|x_1x_2|=\frac\pi2 \text{ for all } x_i\in X_i.$$
Let $X_i^\circ$ denote the interior part of $X_i$
\footnote{At any $p\in X^\circ$ and $\partial X$, $\Sigma_pX$ has empty
and non-empty boundary respectively ([BGP]).}.
For $p_1\in X_1^\circ$ and any $p_2\in X_2$ and any $[p_1p_2]$,
it is easy to see that
$$\uparrow_{p_1}^{p_2}\in (\Sigma_{p_1}X_1)^\perp.\tag{1.1}$$
Since $|x_1p_2|=\frac\pi2$ for all $x_1\in X_1$, by the first
variation formula ([BGP]) it follows that $\uparrow_{p_1}^{p_2}\in
(\Sigma_{p_1}X_1)^{\geq\frac{\pi}{2}}$. Then (1.1) follows from
Lemma 1.3 (note that $p_1\in X_1^\circ$ and $X_1$ is convex, thus
$\Sigma_{p_1}X_1$ is a convex subset with an empty boundary).
Similarly, for any $[p_1p_2]$ with $p_1\in X_1$ and $p_2\in
X_2^\circ$,  we have that
$$\uparrow_{p_2}^{p_1}\in (\Sigma_{p_2}X_2)^\perp.\tag{1.2}$$

By (1.2) and (1.1.3), we have the following property.

\proclaim{Lemma 1.4} For $p_1\in X_1$, $p_2\in X_2^\circ$, given
$[p_1p_2]$ and $[p_2p_2'] \subset X_2$, there exists a $[p_1p_2']$
such that $[p_1p_2]$, $[p_2p_2']$ and $[p_1p_2']$ bound a
convex surface which can be embedded isometrically into the unit sphere
$\Bbb S^2$.
\endproclaim

We make a convention that a convex surface in $X$
as in Lemma 1.4 is called a {\it convex  spherical surface}.
And in the rest, if $\dim(X_i)=0$, we let $X_i^\circ=X_i$.

Now, we assume that $\dim(X_2)>0$, and let $p_2\in
X_2^\circ$ and $p_1\in X_1^\circ$. Let $Q(p_2)$ denote the `cut locus' of $p_2$ in $X_2$,
i.e., $Q(p_2)=\{x\in X,\, [p_2x] \text{ is a maximal minimal geodesic}\}$.
For a fixed $[p_1p_2]$, given $[p_2p_2']$ with $p_2'\in Q(p_2)$, by
(1.2) and Lemma 1.4 there is a convex spherical surface $S\supset [p_1p_2], [p_2p_2']$
(note that the convex spherical surface may not be unique).
Let $c(t)|_{t\in [0,1]}=[p_2p_2']$ with $c(0)=p_2$ and $c(1)=p_2'$, \
and $[p_1c(t)]\subset S$. We define a map
$$f_{[p_1p_2]}(c(t))=\uparrow_{p_1}^{c(t)}\,\, \in (\Sigma_{p_1}X_1)^\perp\ \text{(see (1.1))}.$$
Note that $c(t)$ with $0<t<1$ is in the unique maximal minimal geodesic $[p_2p_2']$,
so $f_{[p_1p_2]}$ is well-defined on $X_2-Q(p_2)$. But it is possible that
there are more than one minimal geodesic between $p_2$ and $p_2'$. Anyway, by choosing
an arbitrary $f_{[p_1p_2]}(p_2')$ for $p_2'\in Q$, we can define a map
$$f_{[p_1p_2]}:X_2\to (\Sigma_{p_1}X_1)^\perp.$$
From the definition, it is clear that for any $p_2'\in X_2$,
$$|f_{[p_1p_2]}(p_2)f_{[p_1p_2]}(p_2')|=|p_2p_2'|,\tag 1.3$$
and that if there is a unique minimal geodesic between $p_2$ and $p_2'$, then
$f_{[p_1p_2]}([p_2p_2'])$ is a minimal geodesic in $(\Sigma_{p_1}X_1)^\perp$.
Consequently, $f_{[p_1p_2]}$ naturally induces a {\it tangent map}
$$\text{\rm D}f_{[p_1p_2]}:\Sigma_{p_2}X_2
\longrightarrow\Sigma_{\uparrow_{p_1}^{p_2}}(\Sigma_{p_1}X_1)^\perp.$$

On the other hand, by (1.1.2) we have that
$$|f_{[p_1p_2]}(p_2')f_{[p_1p_2]}(p_2'')|\geq|p_2'p_2''|\tag 1.4$$
for all $p_2',p_2''\in X_2$.

From the above discussion, we conclude the following.

\proclaim{Proposition 1.5} {\rm (1.5.1)} $f_{[p_1p_2]}$ is
distance non-decreasing; and $f_{[p_1p_2]}$ maps a
$(n_2,\delta)$-strainer at $p_2$ to a $(n_2,\delta)$-strainer at
$\uparrow_{p_1}^{p_2}$, and thus
$\dim(X_2)\leq\dim((\Sigma_{p_1}X_1)^\perp)$.

\vskip1mm

\noindent{\rm (1.5.2)} $\text{\rm D}f_{[p_1p_2]}$ is also a distance
non-decreasing map.
\endproclaim

\demo{Proof} By (1.4), $f_{[p_1p_2]}$ is distance
non-decreasing. Then by (1.3) and (1.1.2), $f_{[p_1p_2]}$ preserves an
$(n_2,\delta)$-strainer (see [BGP] for the definition of the
$(n,\delta)$-strainer), and $\text{D}f_{[p_1p_2]}$ is distance non-decreasing.
\hfill$\qed$\enddemo

The observation that both $f_{[p_1p_2]}$ and $\text{\rm D}f_{[p_1p_2]}$  are
distance nondecreasing maps will play a key rule in the proof of
Theorem A.

\vskip2mm

Based on the geometry of two $\frac \pi2$-part convex subsets in
Lemma 1.2 and 1.4, and properties of $f_{[p_1p_2]}$ and
$\text{D}f_{[p_1p_2]}$ in Proposition 1.5, we will establish the following
key properties that are used in the proof of Theorem A (see the outline
of the proof at the end of Introduction).

Let $\lambda_{x_1x_2}$ be the number of minimal geodesics between $x_1$ and $x_2$.

\proclaim{Key Lemma 1.6} Let $X\in \text{Alex}^n(1)$, and let $X_1$
and $X_2$ be its compact convex subsets with
$|x_1x_2|=\frac\pi2$ for all $x_i\in X_i$ and
$\dim(X_1)+\dim(X_2)=n-1$. Then $\lambda_{p_1p_2}$ has a maximum $m$
for $p_i\in X_i^\circ$; moreover, $X_i^m\triangleq\{p\in X_i^\circ|\
\exists\ q\in X_j^\circ\ (j\neq i)\ \text{s.t. } \lambda_{pq}=m\}$
is open and dense in $X_i^\circ$, and $\lambda_{p_1p_2}=m$ for
all $p_i\in X_i^m$.
\endproclaim

\proclaim{Key Lemma 1.7} Let the assumptions be as in Key Lemma 1.6.
For any $[p_1p_2]$ with $p_i\in X_i^\circ$, we have that
$\uparrow_{p_1}^{p_2}\in ((\Sigma_{p_1}X_1)^\perp)^\circ$. In particular,
if $\partial X_2=\emptyset$, then $\partial(\Sigma_{p_1}X_1)^\perp=\emptyset$,
and $(\Sigma_{p_1}X_1)^\perp=\Uparrow_{p_1}^{X_2}$.
\endproclaim

We conclude this section with an application of
Proposition 1.5, which together with Lemma 1.2 implies (A1).

\proclaim{Theorem 1.8} Let $X\in \text{Alex}^n(1)$, and let $X_1$
and $X_2$ be its two compact convex subsets with
$|x_1x_2|=\frac\pi2$ for all $x_i\in X_i$. Then
$\dim(X_1)+\dim(X_2)\le n-1$.
\endproclaim

\demo{Proof} We observe that Theorem 1.8 holds if $n=1$, and proceed
by induction on $n$. For $p_1\in X_1^\circ$, the following subsets of
$\Sigma_{p_1}X$ are convex:
$$\Sigma_{p_1}X_1, \qquad  (\Sigma_{p_1}X_1)^\perp\ (=(\Sigma_{p_1}X_1)^{=\frac\pi2}).\tag{1.5}$$
(Here, if $\dim(X_1)=0$, then $\Sigma_{p_1}X_1=\emptyset$ and
$\dim(\Sigma_{p_1}X_1)=-1$, and
$(\Sigma_{p_1}X_1)^\perp=\Sigma_{p_1}X$. Note that $\Sigma_{p_1}X_1$
is convex and has an empty boundary because $p_1\in X_1^\circ$ and
$X_1$ is convex, so by Lemma 1.3
$(\Sigma_{p_1}X_1)^{=\frac\pi2}=(\Sigma_{p_1}X_1)^{\geq\frac\pi2}$,
which is convex by (1.1.1).) Since $\Sigma_{p_1}X\in
\text{Alex}^{n-1}(1)$, by induction we have that
$$\dim(\Sigma_{p_1}X_1)+\dim((\Sigma_{p_1}X_1)^\perp)\leq \dim(\Sigma_{p_1}X)-1=n-2,$$ which implies
$$\dim(X_1)+\dim((\Sigma_{p_1}X_1)^\perp)\leq n-1.$$
On the other hand, by (1.5.1) we have that
$$\dim(X_2)\leq\dim((\Sigma_{p_1}X_1)^\perp).$$ Hence, $\dim(X_1)+\dim(X_2)\leq n-1$.
\hfill$\qed$
\enddemo

\remark{Remark \rm1.9} From the above proof, it is clear that
``$\dim(X_1)+\dim(X_2)=n-1$'' implies that
$\dim((\Sigma_{p_i}X_i)^\perp)=\dim(X_j)\ (j\neq i)$ for $p_i\in
X_i^\circ$, and that $(\Sigma_{p_i}X_i)^\perp$ consists of either a
point or two points with distance $\pi$ if
$\dim((\Sigma_{p_i}X_i)^\perp)=0$.
\endremark

%%%%%%%%%%%%%%%%%%%%%%%%%%%%%%%%%%% Section 2 %%%%%%%%%%%%%%%%%%%%%%%%%%%%%%%%%%%%%%%%%%%

\head 2. Proof of Key Lemma 1.6
\endhead

We first point out that Key Lemma 1.6 is obvious if $X_1=\{p_1,p_1'\}$
(or $X_2=\{p_2,p_2'\}$). It is our convention that $|p_1p_1'|=\pi$,
and thus $X=X_1*N$ and $X_2\subseteq N
\in\text{Alex}^{n-1}(1)$. Hence, $\lambda_{p_1p_2}=1$.

In the following proof, we need only to consider two cases:
$X_1$ is a point, and $\dim(X_i)>0$, $i=1, 2$. Let's first
bound $\lambda_{p_1p_2}$ in Key Lemma 1.6.

\proclaim{Lemma 2.1} {\rm (2.1.1)} If $X_1=\{p_1\}$, then
$\lambda_{p_1p_2}\leq 2$ for any $p_2\in X_2^\circ$.

\vskip1mm

\noindent {\rm(2.1.2)} If $\dim(X_i)>0$ for $i=1$ and $2$, then
there exists a positive number $m$ such that $\lambda_{p_1p_2}\leq
m$ for any $p_i\in X_i^\circ$.
\endproclaim

\demo{Proof} (2.1.1) For any $p_2\in X_2^\circ$, by Remark 1.9
$\dim((\Sigma_{p_2}X_2)^\perp)=0$ (note that $\dim(X_2)=n-1$),
and thus $(\Sigma_{p_2}X_2)^\perp$ is either a point or two points with
distance $\pi$. Since $\uparrow_{p_2}^{p_1}\in (\Sigma_{p_2}X_2)^\perp$
for any $[p_1p_2]$ (see (1.2)), we have that $\lambda_{p_1p_2}\leq 2$.

\vskip1mm

(2.1.2) By (2.1.1), we are able to apply the inductive argument starting
with $\dim(X_1)=0$. Let $n_i=\dim(X_i)$.

We first give the proof for a special case where $p_i\in X_i^\circ$ is an
$(n_i,\delta)$-strained point (i.e. a point with an
$(n_i,\delta)$-strainer).

Firstly, under an additional assumption that there is
$(n_i,\delta)$-strained point $q_i\in X_i$ such that
$\lambda_{q_1q_2}<\infty$, we prove that $\lambda_{p_1p_2}$ has a
maximum $m$. If this is not true, then there is a sequence of
$(n_i,\delta)$-strained points $p_i^k\in X_i$ such that
$\lim\limits_{k\to\infty}\lambda_{p_1^kp_2^k}=\infty$. Now we
consider $\lambda_{q_1p_2^k}$. If
$\lim\limits_{k\to\infty}\lambda_{q_1p_2^k}<\infty$ (resp.
$\lim\limits_{k\to\infty}\lambda_{q_1p_2^k}=\infty$) passing a
subsequence, then by Lemma 1.4 there is $[q_1p_2^k]$ and
$[q_1p_1^k]$ (resp. $[q_2q_1]_k$ and $[q_2p_2^k$], where
$[q_2q_1]_k$ is some minimal geodesic between $q_2$ and $q_1$) such
that
$\lim\limits_{k\to\infty}\lambda_{\uparrow_{q_1}^{p_2^k}\uparrow_{q_1}^{p_1^k}}=\infty$
(resp.
$\lim\limits_{k\to\infty}\lambda_{(\uparrow_{q_2}^{q_1})_k\uparrow_{q_2}^{p_2^k}}=\infty$)
in $\Sigma_{q_1}X$ (resp. $\Sigma_{q_2}X$). On the other hand,
according to Remark 1.9, $\dim((\Sigma_{q_1}X_1)^\perp)=n_2$ and
$\dim((\Sigma_{q_2}X_2)^\perp)=n_1$; and by Proposition 1.5,
$\uparrow_{q_1}^{p_2^k}$ and $(\uparrow_{q_2}^{p_1})_k$ are
$(n_2,\delta)$- and $(n_1,\delta)$-strained points in
$(\Sigma_{q_1}X_1)^\perp$ and $(\Sigma_{q_2}X_2)^\perp$
respectively. This implies that $\uparrow_{q_1}^{p_2^k}\in
((\Sigma_{q_1}X_1)^\perp)^\circ$ and $(\uparrow_{q_2}^{q_1})_k\in
((\Sigma_{q_2}X_2)^\perp)^\circ$. Moreover,
$\uparrow_{q_i}^{q_i^k}\in (\Sigma_{q_i}X_i)^\circ$ because
$\Sigma_{q_i}X_i$ has an empty boundary. Then, by the inductive
assumption on $\Sigma_{q_1}X$ (resp. $\Sigma_{q_2}X$), we can
conclude that
$\lambda_{\uparrow_{q_1}^{p_2^k}\uparrow_{q_1}^{p_1^k}}$ (resp.
$\lambda_{(\uparrow_{q_2}^{q_1})_k\uparrow_{q_2}^{p_2^k}}$) have an
upper bound for all $k$, a contradiction.

Secondly, we verify the above additional assumption, i.e., there is
$(n_i,\delta)$-strained point $q_i\in X_i$ such that $\lambda_{q_1q_2}<\infty$.
If not true, then for any $(n_i,\delta)$-strained point $p_i\in X_i$,
$\lambda_{p_1p_2}=\infty$. Note that ``$\lambda_{p_1p_2}=\infty$''
implies that there are minimal geodesics
$[p_1p_2]$ and $\{[p_1p_2]_j\}_{j=1}^\infty$ between $p_1$ and $p_2$
such that
$$\lim_{j\to \infty}[p_1p_2]_j=[p_1p_2]\text{ (so
$\lim_{j\to\infty}(\uparrow_{p_1}^{p_2})_j=\uparrow_{p_1}^{p_2}$)}.
$$
Now we consider the maps $f_{[p_1p_2]}$ and $\text{\rm
D}f_{[p_1p_2]}$ (defined in Section 1). For convenience, we let
$\zeta_j$ and $\zeta$ denote $(\uparrow_{p_1}^{p_2})_j$ and
$\uparrow_{p_1}^{p_2}$ respectively. By Proposition 1.5,
$f_{[p_1p_2]}$ maps an $(n_2,\delta)$-strainer at $p_2$ in $X_2$ to
an $(n_2,\delta)$-strainer at $\zeta$ in $(\Sigma_{p_1}X_1)^\perp$.
On the other hand, we have $\dim((\Sigma_{p_1}X_1)^\perp)=n_2$ (see
Remark 1.9). It then follows that both
$\Sigma_{p_2}X_2$ and $\Sigma_{\zeta}(\Sigma_{p_1}X_1)^\perp$ are
$\chi(\delta)$-isometric to the unit sphere $\Bbb S^{n_2-1}$, where $\chi(\delta)\to
0$ as $\delta\to0$ (ref. [BGP]). Together with that $\text{\rm
D}f_{[p_1p_2]}$ is distance nondecreasing (see Proposition 1.5),
this implies that for any
$\eta\in\Sigma_{\zeta}(\Sigma_{p_1}X_1)^\perp$ there exists
$\xi_\eta\in\Sigma_{p_2}X_2$ such that $|\text{\rm
D}f_{[p_1p_2]}(\xi_\eta)\eta|<\chi(\delta)$. Moreover, we can assume
that there is a $[p_2q_\eta]$ in $X_2$ such that
$\xi_\eta=\uparrow_{p_2}^{q_\eta}$, so $\text{\rm
D}f_{[p_1p_2]}(\xi_\eta)$ is the direction of
$f_{[p_1p_2]}([p_2q_\eta])$ at $\zeta$. Since
$\lim_{j\to\infty}\zeta_j=\zeta$, we assume that
$\uparrow_\zeta^{\zeta_j}$ converges to some
$\bar\eta\in\Sigma_{\zeta}(\Sigma_{p_1}X_1)^\perp$ passing to a
subsequence. Hence, for sufficiently large $j_0$, we can assume that
$|\text{\rm
D}f_{[p_1p_2]}(\xi_{\bar\eta})\uparrow_\zeta^{\zeta_{j_0}}|<2\chi(\delta)$,
and select $\zeta_0\triangleq f_{[p_1p_2]}(p_2^0)$ with
$p_2^0\in[p_2q_{\bar\eta}]$ such that
$|\zeta\zeta_0|=|\zeta\zeta_{j_0}|$ (note that
$|\zeta\zeta_0|=|p_2p_2^0|$). By applying (1.1.2) on the hinge
formed by $[\zeta\zeta_0]$ and $[\zeta\zeta_{j_0}]$, we conclude
that
$$|\zeta_{j_0}\zeta_0|<2\chi(\delta)|\zeta\zeta_0|.$$
However, by applying (1.1.2) on the hinge formed by $[p_1p_2]_{j_0}$
and $[p_1p_2^0]$ (with $\uparrow_{p_1}^{p_2^0}=\zeta_0$), we
conclude that
$$|\zeta_{j_0}\zeta_0|\geq |p_2p_2^0|=|\zeta\zeta_0|,$$
a contradiction. By now, the proof for the special case is complete.

\vskip1mm

Next, we will prove that $\lambda_{p_1p_2}\leq m$ when one of
$p_i$ is $(n_i,\delta)$-strained point, say $p_1$. If $\lambda_{p_1p_2}>m$,
we select minimal geodesics
$\{[p_1p_2]_k\}_{k=1}^{m+1}$, and let
$$\varepsilon=\min\limits_{1\leq i\neq j\leq
m+1}\{|(\uparrow_{p_1}^{p_2})_i(\uparrow_{p_1}^{p_2})_j|\}.$$ Then by Lemma 2.2 below, we
get that $\lambda_{p_1p_2'}\geq m+1$ for any $p_2'\in
B(p_2,\frac{\varepsilon}{2})\cap X_2$. However, there must be an
$(n_2,\delta)$-strained point $p_2'$ in $B(p_2,\frac{\varepsilon}{2})\cap X_2$
(recall that, for any $\nu>0$, there
are $(n_2,\delta)$-strained points in $B(p_2,\nu)\cap X_2$ ([BGP])),
so $\lambda_{p_1p_2'}\leq m$ by the above special case, which contradicts
$\lambda_{p_1p_2'}\geq m+1$.

Finally, for any $p_i\in X_i^\circ$, by repeating the argument in the previous
case we conclude that $\lambda_{p_1p_2}\leq \lambda_{p_1p_2'}\leq m$, where
$p_2'$ is $(n_2,\delta)$-strained point close to $p_2$.

\hfill$\qed$
\enddemo

\proclaim{Lemma 2.2} Let $p_1\in X_1$ and $p_2\in X_2^\circ$,  and
let $[p_1p_2]_1,\cdots,[p_1p_2]_k$ be minimal geodesics between
$p_1$ and $p_2$ with  $\min\limits_{1\leq i\neq j\leq
k}\{|(\uparrow_{p_1}^{p_2})_i(\uparrow_{p_1}^{p_2})_j|\}=\varepsilon$.
Then $\lambda_{p_1p_2'}\geq k$ for any $p_2'\in
B(p_2,\frac{\varepsilon}{2})\cap X_2$.
\endproclaim

\demo{Proof} We will derive a contradiction by assuming that
$\lambda_{p_1p_2'}<k$. By Lemma 1.4, for any $[p_1p_2]_j$ there is a
$[p_1p_2']$ such that
$|(\uparrow_{p_1}^{p_2})_j\uparrow_{p_1}^{p_2'}|=|p_2p_2'|$. Hence,
if $\lambda_{p_1p_2'}<k$, there must be two $[p_1p_2]_j$, say
$[p_1p_2]_1$ and $[p_1p_2]_2$, and one $[p_1p_2']$ such that
$|(\uparrow_{p_1}^{p_2})_j\uparrow_{p_1}^{p_2'}|=|p_2p_2'|\
(j=1,2)$. It then follows that
$$|(\uparrow_{p_1}^{p_2})_1(\uparrow_{p_1}^{p_2})_2|\leq2|p_2p_2'|<\varepsilon,$$
which contradicts ``$\min\limits_{1\leq i\neq j\leq
k}\{|(\uparrow_{p_1}^{p_2})_i(\uparrow_{p_1}^{p_2})_j|\}=\varepsilon$''.
\hfill$\qed$
\enddemo

Based on Lemma 2.1, we will explore the local join structure on $X$,
which will be used to show that $X_i^m$ is open and dense (the latter part of Key Lemma 1.6).

\proclaim{Lemma 2.3} Suppose that the function $\lambda:X_1\times
X_2\to \Bbb N^+$ defined by $\lambda(p,q)=\lambda_{pq}$ attains a
local maximum $m_l$ around $(p_1,p_2)$ with $p_i\in X_i^\circ$. Then
there are neighborhoods $U_i\subset X_i^\circ$ of $p_i$ such that

\noindent{\rm (2.3.1)} $\lambda_{p_1'p_2'}=m_l$ for any $p_1'\in
U_1$ and $p_2'\in U_2$;

\noindent{\rm (2.3.2)} for any
$[p_1p_2]$, both $f_{[p_1p_2]}: U_2\longrightarrow
(\Sigma_{p_1}X_1)^\perp$ and $f_{[p_2p_1]}: U_1\longrightarrow
(\Sigma_{p_2}X_2)^\perp$ are isometrical embedding;

\noindent{\rm (2.3.3)} $U_1*U_2$ can be embedded isometrically into
$X$ around any $[p_1p_2]$.
\endproclaim

\demo{Proof} Let
$$\varepsilon=\min\limits_{1\leq i\neq j\leq
m_l}\{|(\uparrow_{p_2}^{p_1})_i(\uparrow_{p_2}^{p_1})_j|,
|(\uparrow_{p_1}^{p_2})_i(\uparrow_{p_1}^{p_2})_j|\},\tag{2.1}$$
where $\{(\uparrow_{p_2}^{p_1})_i\}$ (resp.
$\{(\uparrow_{p_1}^{p_2})_j\}$) are all directions from $p_2$ to
$p_1$ (resp. from $p_1$ to $p_2$). By Lemma 2.2, for any $p_2'\in
B(p_2,\frac{\varepsilon}{2})\cap X_2$ and $p_1'\in
B(p_1,\frac{\varepsilon}{2})\cap X_1$, we have that
$\lambda_{p_1p_2'}\geq m_l$ and $\lambda_{p_2p_1'}\geq m_l$. On the
other hand, there are neighborhoods $V_i\subset X_i$ of
$p_i$ such that $\lambda|_{V_1\times V_2}\leq m_l$. It then follows
that, for any $p_i'\in V_i\cap B(p_i,\frac{\varepsilon}{2})$,
$$\lambda_{p_1p_2'}=\lambda_{p_1'p_2}=m_l; \tag{2.2}$$
and from the proof of Lemma 2.2 it is easy to
see that, for any fixed $[p_1p_2]$,
$$\exists!\ [p_1p_2'] \text{ s.t. }|\uparrow_{p_1}^{p_2}\uparrow_{p_1}^{p_2'}|=|p_2p_2'|
\text{ and } \ \exists!\ [p_2p_1'] \text{ s.t.
}|\uparrow_{p_2}^{p_1}\uparrow_{p_2}^{p_1'}|=|p_1p_1'|. \tag{2.3}$$
Due to the property of the limit angle (see 2.8.1 in [BGP]),
for sufficiently small convex neighborhood $U_i\subset V_i\cap
B(p_i,\frac{\varepsilon}{2})\cap X_i^\circ$ of $p_i$, we have that
$$\align &\min\limits_{1\leq i\neq j\leq
m_l}\{|(\uparrow_{p_1}^{p_2'})_i(\uparrow_{p_1}^{p_2'})_j|,
|(\uparrow_{p_2'}^{p_1})_i(\uparrow_{p_2'}^{p_1})_j|\}>\frac{\varepsilon}{2}\text{ for any } p_2'\in U_2, \tag{2.4}\\
&\min\limits_{1\leq i\neq j\leq
m_l}\{|(\uparrow_{p_2}^{p_1'})_i(\uparrow_{p_2}^{p_1'})_j|,
|(\uparrow_{p_1'}^{p_2})_i(\uparrow_{p_1'}^{p_2})_j|\}>\frac{\varepsilon}{2}
\text{ for any } p_1'\in U_1.
\endalign$$

\vskip1mm

(2.3.1) Due to (2.4), we can obtain that $\lambda_{p_1'p_2'}=m_l$
(similar to (2.2)) for any $p_i'\in
U_i$ once $U_i$ falls in $B(p_i,\frac{\varepsilon}{4})$.

\vskip2mm

(2.3.2) Since $U_2$ is sufficiently small, we assume that
$U_2\subset B(p_2,\frac{\varepsilon}{8})$. Based on (2.3) and (2.4),
a further observation is that, for any fixed $[p_1p_2]$ and any
$p_2'\in U_2$,
$$\exists!\ [p_1p_2']
\text{ s.t. }|\uparrow_{p_1}^{p_2}\uparrow_{p_1}^{p_2'}|=|p_2p_2'|,
\text{ and } \forall\ p_2^1,p_2^2\in U_2, \
|\uparrow_{p_1}^{p_2^1}\uparrow_{p_1}^{p_2^2}|=|p_2^1p_2^2|.
\tag{2.5}$$ Due to (2.3), we need only to prove that
$|\uparrow_{p_1}^{p_2^1}\uparrow_{p_1}^{p_2^2}|=|p_2^1p_2^2|$. If
this is not true, by Lemma 1.4 there is another
$(\uparrow_{p_1}^{p_2^2})'$ such that
$|\uparrow_{p_1}^{p_2^1}(\uparrow_{p_1}^{p_2^2})'|=|p_2^1p_2^2|$. It
then follows that
$|\uparrow_{p_1}^{p_2^2}(\uparrow_{p_1}^{p_2^2})'|\leq|p_2^2p_2|+
|p_2p_2^1|+|p_2^1p_2^2|<\frac{\varepsilon}2$, which contradicts
(2.4).

Note that (2.5) implies that $f_{[p_1p_2]}: U_2\rightarrow
(\Sigma_{p_1}X_1)^\perp$  is an isometrical embedding, and
similarly, so is $f_{[p_2p_1]}: U_1\rightarrow
(\Sigma_{p_2}X_2)^\perp$ for sufficiently small $U_1$.

\vskip2mm

(2.3.3) We still consider the fixed $[p_1p_2]$, and assume that
$U_i\subset B(p_i,\frac{\varepsilon}{8})$ ($i=1$, 2). For any
$p_2'\in U_2$, let $[p_1p_2']$ be the minimal geodesic determined
in $(2.5)$. Similarly, for any $p_1'\in U_1$,
$$\exists!\ [p_1'p_2']
\text{ s.t.
}|\uparrow_{p_2'}^{p_1}\uparrow_{p_2'}^{p_1'}|=|p_1p_1'|,  \text{
and } \forall\ p_1^1, p_1^2\in U_1,\
|\uparrow_{p_2'}^{p_1^1}\uparrow_{p_2'}^{p_1^2}|=|p_1^1p_1^2|.
\tag{2.6}$$ For the $[p_1p_2]$ and $[p_1'p_2']$,
we define
$${\Bbb D}([p_1p_2],[p_1'p_2'])\triangleq\max\{|xx'|\ |\ x\in[p_1p_2], x'\in[p_1'p_2'], \text{ and } |p_1x|=|p_1'x'|\}.$$
It is not hard to see that $${\Bbb D}([p_1p_2],[p_1'p_2']) \leq
{\Bbb D}([p_1p_2],[p_1p_2'])+{\Bbb
D}([p_1p_2'],[p_1'p_2'])=|p_2p_2'|+|p_1p_1'|.\tag{2.7}$$
Then, due to (2.1), there is $\nu>0$ such that if $|p_2p_2'|+|p_1p_1'|<\nu$, then
$[p_1'p_2']$ is the unique minimal geodesic between $p_1'$ and
$p_2'$ satisfying
$$\Bbb D([p_1p_2],[p_1'p_2'])<\nu.\tag{2.8}$$
Based on (2.8), we claim that {\it for sufficiently small $U_1$
and $U_2$, the $[p_1'p_2^k]$ determined in $(2.6)$ satisfies
$$|\uparrow_{p_1'}^{p_2^1}\uparrow_{p_1'}^{p_2^2}|=|p_2^1p_2^2|.
\tag{2.9}$$ } In fact, by Lemma 1.4 there exists a $[p_1'p_2^1]'$
such that
$|(\uparrow_{p_1'}^{p_2^1})'\uparrow_{p_1'}^{p_2^2}|=|p_2^1p_2^2|$,
so
$$\Bbb D([p_1p_2],[p_1'p_2^1]')\leq \Bbb D([p_1p_2],[p_1'p_2^2])+
\Bbb D([p_1'p_2^2],[p_1'p_2^1]')\leq|p_1p_1'|+|p_2p_2^2|+|p_2^1p_2^2|.$$
Hence, once  $U_1$ and $U_2$ are so small that
$|p_1p_1'|+|p_2p_2^2|+|p_2^1p_2^2|<\nu$, $[p_1'p_2^1]'$ has to be the
$[p_1'p_2^1]$ determined in (2.6). That is, the claim is true.

\vskip1mm

Now, for the $[p_1p_2]$, it suffices to show that
$U_{[p_1p_2]}\triangleq\bigcup_{p_i'\in U_i}[p_1'p_2']$
is isometric to $U_1*U_2$ (where $U_1$ and $U_2$ are in the above claim
and $[p_1'p_2']$ is in (2.6)), i.e.,
for any $x_k\in [p_1^kp_2^k]\subset U_{[p_1p_2]}$ with $p_i^k\in U_i$ ($k=1,2$),
$$|x_1x_2|=|x_1x_2|_*,\tag{2.10}$$
where $|\cdot|_*$ denotes the distance of $U_1*U_2$ (note that $x_k$
can be regarded as the point in the $[p_1^kp_2^k]\subset U_1*U_2$
with $|p_1^kx_k|_*=|p_1^kx_k|$). By (2.9) (and its proof) and
Lemma 1.4, $[p_1^2p_2^1]$, $[p_1^2p_2^2]$ and any $[p_2^1p_2^2]$ bound a
convex spherical surface, denoted by $\overline{\triangle
p_1^2p_2^1p_2^2}$. This implies that $|p_2^1x_2|=|p_2^1x_2|_*$
(similarly, $|p_1^1x_2|=|p_1^1x_2|_*$); and for any
$z\in[p_2^1x_2]\subset\overline{\triangle p_1^2p_2^1p_2^2}$, there
is $p_2^3\in[p_2^1p_2^2]$ such that $z\in[p_1^2p_2^3]\subset
\overline{\triangle p_1^2p_2^1p_2^2}$. Note that
$$\Bbb D([p_1p_2],[p_1^2p_2^3])\leq\Bbb D([p_1p_2],[p_1^2p_2^2])+\Bbb D([p_1^2p_2^2],[p_1^2p_2^3])
\leq |p_1p_1^2|+|p_2p_2^2|+|p_2^2p_2^1|,$$ so it has to hold that
$[p_1^2p_2^3]\subset U_{[p_1p_2]}$ (see (2.8)). Similarly, there is
a convex spherical surface $\overline{\triangle p_1^1p_1^2p_2^3}$
which contains $[p_1^2p_2^3]$, and so $|p_1^1z|=|p_1^1z|_*$ (note
that $z$ can be regarded as the point in the $[p_2^1x_2]\subset
U_1*U_2$ with $|p_2^1z|_*=|p_2^1z|$). Hence, by (1.1.3),
$[p_1^1z]\ (\subset\overline{\triangle p_1^1p_1^2p_2^3})$ and
$[p_2^1x_2]$ determine a convex spherical surface
$\overline{\triangle p_1^1p_2^1x_2}$ ($\supset[p_1^1z]$) bounded by
$[p_2^1x_2]$ and some $[p_1^1p_2^1]'$ and $[p_1^1x_2]$. Similarly,
for any $x\in [p_1^1p_2^1]'$, we have that
$$|x_2x|=|x_2x|_*,$$
which will imply (2.10) as long as we prove that $[p_1^1p_2^1]'$ is just the
$[p_1^1p_2^1]$.  Note that
$[p_1^1z]\subseteq\overline{\triangle p_1^1p_1^2p_2^3}\cap\overline{\triangle p_1^1p_2^1x_2}$
and $[p_2^1z]\subseteq\overline{\triangle p_1^2p_2^1p_2^2}\cap\overline{\triangle p_1^1p_2^1x_2}$,
so it is not hard to see that $\Bbb D([p_1p_2]',[p_1^2p_2^3])<\nu/2$ for sufficiently small $U_1$ and $U_2$
(where $\nu$ is the number in (2.8)). And we can assume that
$\Bbb D([p_1p_2],[p_1^2p_2^3])<\nu/2$ (see (2.7-8)), so
$$\Bbb D([p_1p_2],[p_1^1p_2^1]')\leq \Bbb D([p_1p_2],[p_1^2p_2^3])+
\Bbb D([p_1^2p_2^3],[p_1^1p_2^1]')<\nu.$$ Then, by (2.8), $[p_1^1p_2^1]'$ has to be
the $[p_1^1p_2^1]$ determined in $(2.6)$.
\hfill$\qed$
\enddemo

\remark{Remark \rm2.4} The above argument can be used to prove a
weak version of `local' join structure: for any $[p_1p_2]$ with
$p_i\in X_i^\circ$ and $\lambda_{p_1p_2}<m$, given $[p_1'p_2']$ with
$p_i'\in X_i^\circ$ and $\Bbb D([p_1p_2],[p_1'p_2'])$ sufficiently
small, there is $[p_ip_i']\subset X_i$ such that $[p_1p_2]\cup
[p_1'p_2']$ can be embedded isometrically into
$[p_1p_1']*[p_2p_2']$, i.e., for any $x\in[p_1p_2]^\circ$ and
$x'\in[p_1'p_2']^\circ$,
$$|xx'|=|xx'|_*\quad\text{(similar to (2.10))},$$
where $|\cdot|_*$ denotes the distance of
$[p_1p_1']*[p_2p_2']$ (note that here $[p_1p_1']$ and
$[p_2p_2']$ can be chosen arbitrarily). In what follows,
we point out that, for any $[xx']$, the above
underlying join structure implies a canonical way to
locate the projection of $[xx']$ on $X_i$ \footnote{For any $x\in
[x_1x_2]^\circ$ with $x_i\in X_i$, we say that $x_i$ is {\it the
projection of $x$ on $X_i$}.}, and the projection is
a minimal geodesic $[p_ip_i']$.
Indeed, the way to locate the projection here is the same as to
do this in $[p_1p_1']*[p_2p_2']$ (see below for a step-by-step picture).
\vskip.2in
\epsfxsize=5truein                      % <--- sets horizontal size
\centerline{\epsfbox{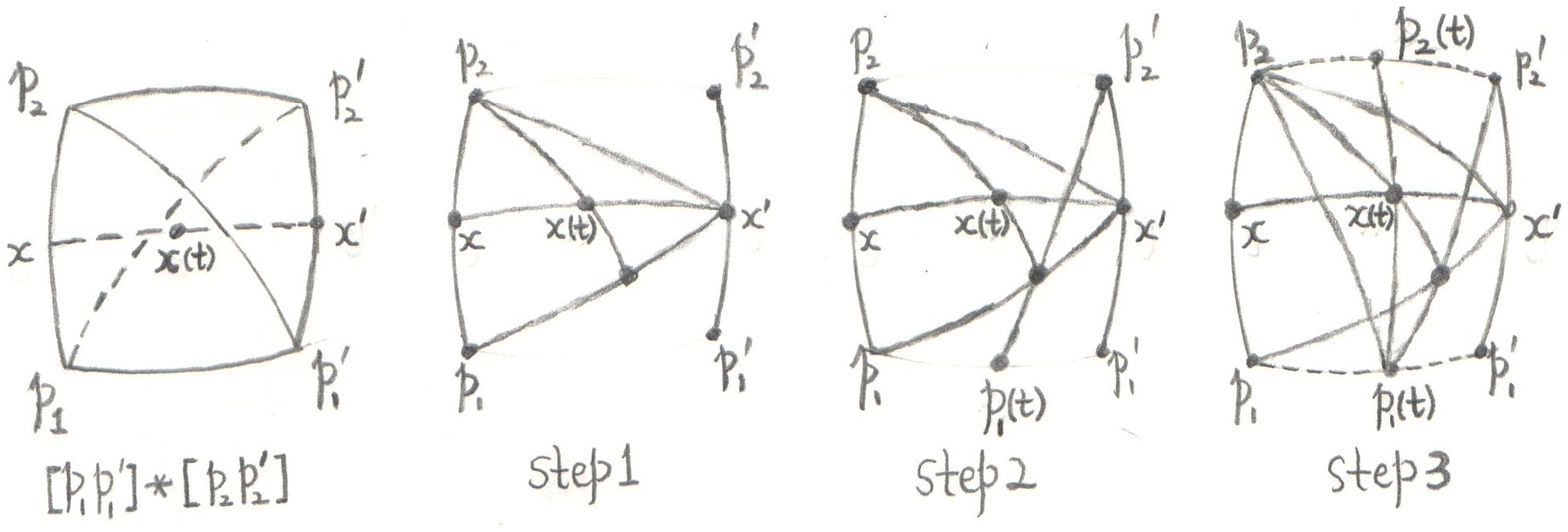}} % <--- centers your figure
\endremark

\demo{Proof of Key Lemma 1.6}

By Lemma 2.1, $\lambda_{p_1p_2}$ has a maximum $m$ for all
$p_i\in X_i^\circ$. And by (2.3.1),  $X_i^m\ (=\{p\in X_i^\circ|\
\exists\ q\in X_j^\circ\ (j\neq i)\ \text{s.t. } \lambda_{pq}=m\})$ is open
in $X^\circ$. Hence, it
suffices to show that $X_i^m$ is dense in $X_i^\circ$ and
that $\lambda_{p_1p_2}=m$ for all $p_i\in X_i^m$. We will give a proof by
induction on $n$, starting with the trivial case, $n=1$.

\vskip1mm

We first prove that for all $p_i\in X_i^m$, $\lambda_{p_1p_2}=m$.
By the definition of $X_i^m$, there is a $p_2'\in X_2^\circ$ and $p_1'\in X_1^\circ$
such that $\lambda_{p_1p_2'}=m=\lambda_{p_1'p_2}$. According to
(2.3.3), there are neighborhoods $U_i\subset X_i^\circ$ of $p_1$ and
$p_2'$ such that $U_1*U_2$ can be embedded isometrically into $X$
around any given $[p_1p_2']$.

We consider $\Sigma_{p_1}X\in\text{Alex}^{n-1}(1)$, in which both
$\Sigma_{p_1}X_1$ and $(\Sigma_{p_1}X_1)^\perp$ are convex (see (1.5)). By
Remark 1.9,
$\dim(\Sigma_{p_1}X_1)+\dim((\Sigma_{p_1}X_1)^\perp)=\dim(\Sigma_{p_1}X)-1$.
Hence, by the inductive assumption on $\Sigma_{p_1}X$,
$\lambda_{\xi\zeta}$ has a maximum $m'$ for any
$\xi\in (\Sigma_{p_1}X_1)^\circ$ ($=\Sigma_{p_1}X_1$)
and $\zeta\in ((\Sigma_{p_1}X_1)^\perp)^\circ$, and
$\lambda_{\xi\zeta}=m'$ for any $\xi\in(\Sigma_{p_1}X_1)^{m'}$
and $\zeta\in((\Sigma_{p_1}X_1)^\perp)^{m'}$. On the other hand,
by (2.3.2), $f_{[p_1p_2']}: U_2\longrightarrow
(\Sigma_{p_1}X_1)^\perp$ is an isometrical embedding, which implies
that $\uparrow_{p_1}^{p_2'}\in ((\Sigma_{p_1}X_1)^\perp)^\circ$
(note that $\dim((\Sigma_{p_1}X_1)^\perp)=\dim(X_2)$ and $p_2'\in
X_2^\circ$). Then that $U_1*U_2$ can be
embedded isometrically into $X$ around $[p_1p_2']$ implies that
$$m'=1.\tag{2.11}$$

Let $p_2''$ be an $(n_2,\delta)$-strained point in $X_2$ closed to $p_2$.
We first show that $\lambda_{p_1p_2''}=m$.
From the proof of (2.1.2), we know that $\uparrow_{p_1}^{p_2''}\in
((\Sigma_{p_1}X_1)^\perp)^\circ$ for any $[p_1p_2'']$ (however, so
far we do NOT know that $\uparrow_{p_1}^{p_2}\in
((\Sigma_{p_1}X_1)^\perp)^\circ$ for a $[p_1p_2]$).
Note that $\lambda_{p_1'p_2''}=m$ by (2.3.1). If $\lambda_{p_1p_2''}<m$,
then by Lemma 1.4
there is some $\uparrow_{p_1}^{p_1'}$ and $\uparrow_{p_1}^{p_2''}$
such that $\lambda_{\uparrow_{p_1}^{p_1'}\uparrow_{p_1}^{p_2''}}>1$,
which contradicts (2.11).

Now, we show that $\lambda_{p_1p_2}=m$.
Since $\lambda_{p_1p_2''}=m$, if $\lambda_{p_1p_2}<m$, then by Lemma
1.4 there are $[p_1p_2'']_k$ ($k=1,2$) and $[p_1p_2]$ such that
$|\uparrow_{p_1}^{p_2}(\uparrow_{p_1}^{p_2''})_k|=|p_2p_2''|$.
Note that there is a unique minimal geodesic between
$(\uparrow_{p_1}^{p_2''})_k\in ((\Sigma_{p_1}X_1)^\perp)^\circ$ and
given $\uparrow_{p_1}^{p_1'}\in\Sigma_{p_1}X$ (see (2.11)),
so by applying Lemma 1.4 on $\Sigma_{p_1}X$ we get that
$$\Bbb D([\uparrow_{p_1}^{p_1'}(\uparrow_{p_1}^{p_2''})_1],
[\uparrow_{p_1}^{p_1'}(\uparrow_{p_1}^{p_2''})_2])=
|(\uparrow_{p_1}^{p_2''})_1(\uparrow_{p_1}^{p_2''})_2|\leq
2|p_2p_2''|.\tag{2.12}$$ On the other hand, still by Lemma 1.4,
there is $[p_1'p_2'']_k$ such that
$[\uparrow_{p_1}^{p_1'}(\uparrow_{p_1}^{p_2''})_k]$
$(\subset \Sigma_{p_1}X)$ is realized by the convex spherical surface $S_k$ bounded
by $[p_1p_2'']_k, [p_1'p_2'']_k$ and $[p_1p_1']$. Let
$x_k\in [p_1'p_2'']_k$ with $|p_1'x_1|=|p_1'x_2|$, and $[p_1x_k]\subset S_k$.
Note that $|p_1x_1|=|p_1x_2|$, so
$|\uparrow_{p_1}^{x_1}\uparrow_{p_1}^{x_2}|\geq|x_1x_2|$ by (1.1.2).
This implies that
$$\Bbb D([\uparrow_{p_1}^{p_1'}(\uparrow_{p_1}^{p_2''})_1],
[\uparrow_{p_1}^{p_1'}(\uparrow_{p_1}^{p_2''})_2]) \geq
\Bbb D([p_1'p_2'']_1,[p_1'p_2'']_2),$$
which contradicts $(2.12)$ (because $|p_2p_2''|$ can be sufficiently small,
but $\Bbb D([p_1'p_2'']_1,$ $[p_1'p_2'']_2)$
has a positive lower bound (see (2.4))).

\vskip1mm

We then prove that $X_i^m$ is dense in  in $X_i^\circ$.
If this is not true, then there is $q_i\in X_i^\circ$ such that
$\lambda_{q_1q_2}=m'<m$ is a local maximum. Similar to (2.11),
we can conclude that there is a unique minimal geodesic between any point in
$\Sigma_{q_i}X_i$ and in $((\Sigma_{q_i}X_i)^\perp)^\circ$.
Now, consider a pair of $p_i\in X_i^m$. By (2.3.1),
$p_i, q_i$ can be chosen to be $(n_i,\delta)$-stained points.
Note that, when $\lambda_{p_1q_2}=m$ (resp. $\lambda_{p_1q_2}<m$),
by Lemma 1.4, there is some $\uparrow_{q_1}^{p_1}\in
\Sigma_{q_1}X_1$ and $\uparrow_{q_1}^{q_2}\in
((\Sigma_{q_1}X_1)^\perp)^\circ$ (resp. $\uparrow_{q_2}^{p_2}\in
\Sigma_{q_2}X_2$ and $\uparrow_{q_2}^{p_1}\in
((\Sigma_{q_2}X_2)^\perp)^\circ$) between which there are at least
two minimal geodesics; a contradiction. \hfill$\qed$
\enddemo

\vskip2mm

We will end this section with some properties of $X_i^m$, which will be used in Section 4.

\proclaim{Corollary 2.5} For any $p_i\in X_i^m$, if
$(\Sigma_{p_i}X_i)^\perp$ has an empty boundary, then
$\Sigma_{p_i}X=\Sigma_{p_i}X_i*(\Sigma_{p_i}X_i)^\perp$.
\endproclaim

Note that Corollary 2.5 can be seen from the proof around (2.11).

\proclaim{Lemma 2.6} For any $p_1\in X_1^\circ$ and
$[p_2p_2']\subset X_2^\circ$, there exists a finite number $m_c$ such that
$\lambda_{p_1p_2''}=m_c$ for all $p_2''\in [p_2p_2']^\circ$, and
$\lambda_{p_1p_2}, \lambda_{p_1p_2'}\leq m_c$.
\endproclaim

\demo{Proof} By Lemma 1.4 (resp. (1.1.3)), any $[p_1p_2]$ (resp. $[p_1p_2'']$)
has to lie in a convex spherical surface bounded by $[p_1p_2]$,
$[p_2p_2']$ and some $[p_1p_2']$. Moreover, the
interior parts of two such surfaces do not intersect. This together
with Lemma 2.1 implies the lemma. \hfill$\qed$
\enddemo

Lemma 2.6 has an immediate corollary.

\proclaim{Corollary 2.7} {\rm (2.7.1)} For any $p_i\in X_i^m$ and
$[p_ip_i']\subset X_i$, $[p_ip_i']\setminus\{p_i'\}$ belongs to
$X_i^m$.

\vskip1mm

\noindent {\rm (2.7.2)} $X_i^m$ is totally convex in $X_i$.
\endproclaim

\head 3. Proof of Key Lemma 1.7
\endhead

In the proof Key Lemma 1.7, we will show that the following multi-valued map
$$f_{p_1}: X_2\to (\Sigma_{p_1}X_1)^\perp\quad \text{by } p_2\mapsto\Uparrow_{p_1}^{p_2}$$
is onto if $\partial X_2=\emptyset$, where $p_1\in X_1^\circ$.
It suffices to check that $f_{p_1}$ is an open map (it is clear that
$f_{p_1}$ is a closed map). Note that at any $(n_2,\delta)$-strained point,
we have already known that $f_{p_1}$ is open (see the arguments after (2.11)).

To prove the openness for $f_{p_1}$, we introduce the concept of {\it cone-neighborhood
isometrical multi-valued map}.

\example{Definition 3.1} Let $X, \tilde X\in \text{Alex}(1)$ with $\dim(\tilde X)\ge\dim(X)\geq1$, and let $f:X\rightarrow \tilde X$ be a
multi-valued map. We say that $f$ is a
{\it cone-neighborhood isometry}  if the following
hold:

\noindent(3.1.1) there is an $m$ such that $\#\{f(p)\}\leq m$ for all
$p\in X$, and that $X_m\triangleq\{p\in X|\ \#\{f(p)\}=m\}$ is dense in $X$.

\noindent(3.1.2) $|\tilde p\tilde q|\geq|pq|$ for any $p, q\in X$,
$\tilde p\in f(p)$ and $\tilde q\in f(q)$; and given any $[pq]$ and
$\tilde p\in f(p)$, there exists a $[\tilde p\tilde q]$
with $\tilde q\in f(q)$ such that $[\tilde p\tilde q]\subseteq f([pq])$ and
$f|_{[pq]}:[pq]\rightarrow [\tilde p\tilde q]$ is an isometry.
\endexample

\remark{Remark \rm3.2} For any $p\in X$ and $\tilde p \in f(p)$ in
Definition 3.1, let $r_{\tilde p}\triangleq\frac{1}{4}\min\{|\tilde
p\tilde p'|\ |\ \tilde p'\in f(p), \tilde p'\neq\tilde p\}$. Then it
follows from (3.1.1) and (3.1.2) that:

\noindent(3.2.1) if $p\in X_m$, then
$f:B(p,r_{\tilde p})\rightarrow B(\tilde p,r_{\tilde p})$ is an isometrical embedding;

\noindent(3.2.2) if $p\in X\setminus X_m$, then
$|\tilde p\tilde q|=|pq|$ for any $q\in B(p,r_{\tilde p})$
and $\tilde q\in f(q)\cap B(\tilde p,r_{\tilde p})$.
\endremark

\vskip2mm

Why do we call such $f$ a cone-neighborhood isometry?
To give an explanation, we first define the {\it $\varepsilon$-cone
neighborhood $V_{[pq],\varepsilon}$} of any given $[pq]\subset X$ as follows:
$$V_{[pq],\varepsilon}\triangleq\left\{x\in X\ |\ |xp|\leq |qp|
\text{ and there is a } [px]  \text{ s.t. } |\uparrow_p^x\uparrow_p^q|<\epsilon\right\}.$$
From Definition 3.1 (and Remark 3.2), for any $\tilde p\in f(p)$ and  $[pq]\subset B(p, r_{\tilde p})$,
it is easy to see that there are minimal geodesics $\{[\tilde p\tilde q_i]\}_{i=1}^k$ ($k\leq m$)
with $\tilde q_i\in f(q)\cap B(\tilde p,r_{\tilde p})$ such that
$$B(\tilde p,r_{\tilde p})\cap f([pq])=\bigcup_{i=1}^k[\tilde p\tilde q_i].\tag{3.1}$$
We know that any triangle $\triangle\tilde x\tilde p\tilde y$, where $\tilde x\in [\tilde p\tilde q_i]$
and $\tilde y\in [\tilde p\tilde q_j]$, will be more and more isometric to its comparison triangle
as $\tilde x,\tilde y\to \tilde p$ ([BGP]).
Hence, it is not hard to see that there is an $N>4$ such that for each $\uparrow_{\tilde p}^{\tilde q_i}$
$$r_{\uparrow_{\tilde p}^{\tilde q_i}}\triangleq\frac{1}{N}\min_{j\neq i}\{|\uparrow_{\tilde p}^{\tilde q_i}\uparrow_{\tilde p}^{\tilde q_j}|\}\tag{3.2}$$
satisfies that, for $\tilde q_{i,\delta}\in [\tilde p\tilde q_i]$ and $q_\delta\in [pq]$ with $|\tilde p\tilde q_{i,\delta}|=|pq_\delta|=\delta$ sufficiently small,
$$V_{[\tilde p\tilde q_{i,\delta}],r_{\uparrow_{\tilde p}^{\tilde q_i}}}\subset \bigcup\limits_{\tilde x\in[\tilde p\tilde q_{i,\delta}]\setminus\{\tilde p\}}B(\tilde x,r_{\tilde x})\tag{3.3}$$
and
$$V_{[pq_\delta],r_{\uparrow_{\tilde p}^{\tilde q_i}}}\subset \bigcup\limits_{x\in[pq_\delta]\setminus\{p\}, \tilde x\in f(x)\cap[\tilde p\tilde q_i]}B(x,r_{\tilde x}).\tag{3.4}$$
It then follows that $$f(V_{[pq_\delta],r_{\uparrow_{\tilde p}^{\tilde q_i}}})
\subseteq V_{[\tilde p\tilde q_{i,\delta}],r_{\uparrow_{\tilde p}^{\tilde q_i}}},$$
due to which we call the multi-valued $f$ a cone-neighborhood isometry.

By the way, we would like to point out that (3.1) implies that
$$[pq]\setminus\{p\}\subset X_m \tag{3.5}$$
for any $[pq]$ with $q\in B(p,r_p)\cap X_m$, where $r_p=\min\{r_{\tilde p}|\tilde p\in f(p)\}$.

\vskip2mm

A substantial property of a cone-neighborhood isometry is:

\proclaim{Proposition 3.3} Let $f:X\rightarrow \tilde X$ be a
cone-neighborhood isometry. Then at any $p\in X$ and
$\tilde p\in f(p)$, $f$ induces a tangent map
$\text{\rm D}f:\Sigma_pX\rightarrow\Sigma_{\tilde p}\tilde X$ such
that $\text{\rm D}f$ is again a cone-neighborhood isometry.
\endproclaim

\demo{Proof} By (3.2.1), if $p\in X_m$, then $f$ induces naturally an isometrical embedding
$\text{\rm D}f:\Sigma_pX\rightarrow\Sigma_{\tilde p}\tilde X$.
Hence, in the rest of the proof, we need only to consider $p\in X\setminus X_m$.

At first, we define $\text{\rm D}f(\eta)$ for all $\eta\in (\Sigma_pX)'$. Recall that $\eta\in (\Sigma_pX)'$
means that there is a $[pq]$ with $q\in B(p,r_{\tilde p})$ such that $\eta=\uparrow_p^q$ ([BGP]).
Since $B(\tilde p,r_{\tilde p})\cap f([pq])=\cup_{i=1}^k[\tilde p\tilde q_i]$ ($k\leq m$) where $\tilde q_i\in f(q)\cap B(\tilde p,r_{\tilde p})$ (see (3.1)),
we define $$\text{\rm D}f(\eta)=\{\uparrow_{\tilde p}^{\tilde q_i}\}_{i=1}^k,\tag{3.6}$$
on which we will first give two claims before defining $\text{\rm D}f(\zeta)$ for $\zeta\not\in (\Sigma_pX)'$.

{\bf Claim 1}: {\it there is an $m'$ such that $\#\{\text{\rm D}f(\eta)\}\leq m'$ for all $\eta\in (\Sigma_pX)'$;
and if there is a $[pq]$ with $q\in B(p,r_{\tilde p})\cap X_m$ such that $\eta=\uparrow_p^q$, then $\#\{\text{\rm D}f(\eta)\}=m'$.}
In order to see the claim, we fix a point $q\in B(p,r_{\tilde p})\cap X_m$,
and notice that $\#\{\text{\rm D}f(\eta)\}=\#\{f(q)\cap B(\tilde p,r_{\tilde p})\}$.
By (3.1.2) and the choice of $r_{\tilde p}$ (in Remark 3.2),
for any $[zq]$ with $z\in B(p,r_{\tilde p})$ and $\tilde z\in f(z)\cap B(\tilde p,r_{\tilde p})$,
there exists a $[\tilde z\tilde q_i]$ with $\tilde q_i\in f(q)\cap B(\tilde p,r_{\tilde p})$
such that $f|_{[zq]}:[zq]\rightarrow [\tilde z\tilde q_i]$ is an isometry.
It follows that $\#\{f(z)\cap B(\tilde p,r_{\tilde p})\}\leq\#\{f(q)\cap B(\tilde p,r_{\tilde p})\}$
(otherwise, for some $\tilde q_i\in f(q)\cap B(\tilde p,r_{\tilde p})$, $f([zq])\cap B(\tilde q_i,r_{\tilde q_i})$ consists of
at least two geodesics starting from $\tilde q_i$, which contradicts (3.2.1) because $q\in X_m$).
Note that this also implies that $\#\{f(z)\cap B(\tilde p,r_{\tilde p})\}=\#\{f(q)\cap B(\tilde p,r_{\tilde p})\}$
if $z$ also belongs to $X_m$. Therefore, Claim 1 follows.

{\bf Claim 2}: {\it Let $q,z\in B(p,r_{\tilde p})$. For any  $[pq]$, $[pz]$ and
$[\tilde p\tilde z]\subset f([pz])\cap B(\tilde p,r_{\tilde p})$, there is a $[\tilde p\tilde q_i]\subset f([pq])\cap B(\tilde p,r_{\tilde p})$
such that $|\uparrow_{\tilde p}^{\tilde z}\uparrow_{\tilde p}^{\tilde q_i}|=|\uparrow_{p}^z\uparrow_{p}^q|$.}
Let $x_j\in [pq]$ and $y_j\in [pz]$ with $x_j,y_j\to p$ as $j\to\infty$.
We select $\tilde y_j\in f(y_j)\cap [\tilde p\tilde z]$.
Similarly, for any $[y_jx_j]$, there is a $[\tilde y_j\tilde x_j]$ with $\tilde x_j\in f(x_j)\cap B(\tilde p, r_{\tilde p})$
such that $f|_{[y_jx_j]}:[y_jx_j]\rightarrow [\tilde y_j\tilde x_j]$ is an isometry.
Passing to a subsequence, we can assume that $\tilde x_j$ falls in some
$[\tilde p\tilde q_i]\subset f([pq])\cap B(\tilde p,r_{\tilde p})$.
It then follows that
$|\uparrow_{\tilde p}^{\tilde z}\uparrow_{\tilde p}^{\tilde q_i}|=|\uparrow_{p}^z\uparrow_{p}^q|$
(note that $|\tilde p\tilde x_j|=|px_j|, |\tilde p\tilde y_j|=|py_j|$ and $|\tilde x_j\tilde y_j|=|x_jy_j|$).

\vskip1mm

Next, we will define $\text{\rm D}f(\zeta)$ for all $\zeta\in
\Sigma_pX\setminus(\Sigma_pX)'$. Note that we can
select $[pz_j]$ with $z_j\in B(p, r_{\tilde p})\cap X_m$ such that
$\uparrow_p^{z_j}\to \zeta$ as $j\to\infty$.

{\bf Claim 3}: {\it $B\triangleq\cup_{j=1}^\infty \text{\rm
D}f(\uparrow_p^{z_j})$ has at most $m'$ limiting points in
$\Sigma_{\tilde p}\tilde X$, where $m'$ is the number in Claim 1;
and the limiting points of $B$ do not depend on the choice of
$[pz_j]$.} If the claim is not true, then we can let $\tilde
\zeta_1, \cdots, \tilde \zeta_{m'+1}$ be the limiting points of $B$.
Note that for any $0<\epsilon\ll\min_{1\leq i\neq i'\leq
m'+1}\{|\tilde \zeta_i\tilde \zeta_{i'}|\}$, there is a $J>0$ such
that $|\uparrow_p^{z_j}\uparrow_p^{z_{j'}}|\leq\epsilon$ if
$j,j'>J$. This together with Claim 2 implies that, for any fixed
$j>J$ and $\tilde \zeta_i$ ($1\leq i\leq m'+1$), there is a
$\tilde\xi\in \text{\rm D}f(\uparrow_p^{z_j})$ such that
$|\tilde\xi\tilde\zeta_i|\leq\epsilon$. I.e., it has to hold that
$\#\{\text{\rm D}f(\uparrow_p^{z_j})\}>m'$, which contradicts Claim
1. Hence, $B$ has at most $m'$ limiting points; moreover, we can
similarly conclude that the limiting points of $B$ do not depend on
the choice of $[pz_j]$.

Based on Claim 3, for $\zeta\in\Sigma_pX\setminus(\Sigma_pX)'$,
we define $$\text{\rm D}f(\zeta)=\{\text{the limiting points of } B\}.$$
So far, we have finished the definition of a multi-valued map $\text{\rm D}f:\Sigma_pX\to\Sigma_{\tilde p}\tilde
X$.

\vskip1mm

At last, we need only to check that $\text{\rm D}f$ satisfies the corresponding
(3.1.1-2).

\vskip0.5mm

``(3.1.1)'': By Claim 1 and 3, it is clear that $\#\{\text{\rm D}f(\eta)\}\leq m'$ for all
$\eta\in \Sigma_pX$; and by the latter part of Claim 1, $(\Sigma_pX)_{m'}\triangleq\{\eta\in \Sigma_pX|\ \#\{\text{\rm D}f(\eta)\}=m'\}$
is dense in $\Sigma_pX$ because $X_m$ is dense in $X$ and $(\Sigma_pX)'$ is dense in $\Sigma_pX$ ([BGP]).

\vskip0.5mm

``(3.1.2)'': For any $\eta,\zeta\in \Sigma_pX$, we need to prove that
$|\tilde \eta\tilde \zeta|\geq|\eta\zeta|$ for any $\tilde \eta\in \text{\rm D}f(\eta)$ and  $\tilde \zeta\in \text{\rm D}f(\zeta)$, and that
for any given $[\eta\zeta]$ and $\tilde \eta_0\in \text{\rm D}f(\eta)$ there is a $[\tilde\eta_0\tilde\zeta_0]$
with $\tilde \zeta_0\in \text{\rm D}f(\zeta)$ such that $\text{\rm D}f|_{[\eta\zeta]}:[\eta\zeta]\to[\tilde\eta_0\tilde\zeta_0]$ is an isometry.

We first consider $\eta,\zeta\in (\Sigma_pX)'$, i.e., there are $[pq],[pz]\subset B(p,r_{\tilde p})$ such that
$\eta=\uparrow_p^q$ and $\zeta=\uparrow_p^z$.
From (3.6), there is a $[\tilde p\tilde q]$ and $[\tilde p\tilde z]$ with $\tilde q\in f(q)\cap B(\tilde p,r_{\tilde p})$
and $\tilde z\in f(z)\cap B(\tilde p,r_{\tilde p})$ such that $\tilde \eta=\uparrow_{\tilde p}^{\tilde q}$
and $\tilde \zeta=\uparrow_{\tilde p}^{\tilde z}$. Note that for any
$\tilde q'\in[\tilde p\tilde q]$ and $q'\in f^{-1}(\tilde q')$ and any $\tilde z'\in[\tilde p\tilde z]$
and $z'\in f^{-1}(\tilde z')$, we have $|\tilde q'\tilde z'|\geq|q'z'|$ (see (3.1.2)),
$|\tilde p\tilde q'|=|pq'|$ and $|\tilde p\tilde z'|=|pz'|$.
It follows that $|\uparrow_{\tilde p}^{\tilde q}\uparrow_{\tilde p}^{\tilde z}|\geq|\uparrow_p^q\uparrow_p^z|$, i.e.,
$$|\tilde \eta\tilde \zeta|\geq|\eta\zeta|.$$
Next we will find a $[\tilde\eta_0\tilde\zeta_0]$
such that $\text{\rm D}f|_{[\eta\zeta]}:[\eta\zeta]\to[\tilde\eta_0\tilde\zeta_0]$ is an isometry.
Note that $[\eta\zeta]$ can be determined by a sequence of triangles $\{\triangle pq_jz_j\}_{j=1}^\infty$
with $q_j\in [pq]$ and $z_j\in[pz]$ and $|pq_j|=|pz_j|\to0$ as $j\to\infty$
(note that $\triangle pq_jz_j$ in $(\frac{1}{|pq_j|^2}X,p)$ converges to the cone $C([\eta\zeta])\subset T_p=C(\Sigma_pX)$ ([BGP])).
Now, we select $[\tilde p\tilde q_0]$ with $\tilde q_0\in f(q)\cap B(\tilde p,r_{\tilde p})$ such that $\tilde \eta_0=\uparrow_{\tilde p}^{\tilde q_0}$.
By (3.1.2), for $\tilde q_j\triangleq f(q_j)\cap[\tilde p\tilde q_0]$, there is a $[\tilde q_j\tilde z_j]$ with $\tilde z_j\in f(z_j)$
such that $f|_{[q_jz_j]}:[q_jz_j]\ (\subset\triangle pq_jz_j)\to [\tilde q_j\tilde z_j]$ is an isometry.
Since $f([pz])\cap B(\tilde p,r_{\tilde p})$ consists of at most $m'$ pieces of geodesics starting from $\tilde p$ (see Claim 1),
passing to a subsequence we can assume that all $\tilde z_j$ lie in a $[\tilde p\tilde z_0]$ with $\tilde z_0\in f(z)\cap B(\tilde p,r_{\tilde p})$.
It is not hard to see that a subsequence of triangles $\{\triangle \tilde p\tilde q_j\tilde z_j\}_{j=1}^\infty$ (where $\triangle \tilde p\tilde q_j\tilde z_j$ is formed by
$[\tilde p\tilde q_j]$, $[\tilde p\tilde z_j]$ and $[\tilde q_j\tilde z_j]$) determines a $[\tilde\eta_0\tilde\zeta_0]$
with $\zeta_0=\uparrow_{\tilde p}^{\tilde z_0}$ in $\Sigma_{\tilde p}\tilde X$
such that $\text{\rm D}f|_{[\eta\zeta]}:[\eta\zeta]\to[\tilde\eta_0\tilde\zeta_0]$ is an isometry.

Now we assume that at least one of $\eta$ and $\zeta$ does not belong to $(\Sigma_pX)'$, say $\eta,\zeta\not\in(\Sigma_pX)'$.
Since we can find $\eta_j,\zeta_j\in (\Sigma_pX)'$
such that $\eta_j\to \eta$ and $\zeta_j\to \zeta$ as $j\to \infty$, we can use the standard limiting argument (together with
what we have proven in the above) to complete the proof.
\hfill$\qed$
\enddemo

Similar to a (single-valued) continuous map, any cone-neighborhood
isometry maps a compact set to a compact set.

\proclaim{Proposition 3.4} Let $f:X\rightarrow \tilde X$ be a
cone-neighborhood isometry. If $X$ is compact, then $f(X)$
is also compact (and thus $f(X)$ is closed in $\tilde X$).
\endproclaim

\demo{Proof} It suffices to show that any sequence $\{\tilde x_i\}_{i=1}^\infty\subset f(X)$
contains a subsequence which converges to a point in $f(X)$
(that $f(X)$ is closed in $\tilde X$ is because $\tilde X$ is Hausdorff and $f(X)$
is compact). Let $x_i=f^{-1}(\tilde x_i)$ for all $i$.
Since $X$ is compact, passing to a subsequence we can assume that $x_i\to x$ as $i\to\infty$.
Then due to (3.1.2), $\{\tilde x_i\}_{i=1}^\infty$ has to contain a subsequence which
converges to some point of $f(x)$ (note that $\#\{f(x)\}\leq m$).
\hfill$\qed$
\enddemo

Based on Proposition 3.3 and 3.4, we can derive the following important property of
a cone-neighborhood isometry.

\proclaim{Theorem 3.5} Let $X, \tilde X\in \text{Alex}^n(1)$ with $n\geq1$, and let $f:X\rightarrow \tilde X$ be a
cone-neighborhood isometry. If $X$ is compact and $\partial X=\emptyset$,
then $\tilde X=f(X)$ (which is compact by Proposition 3.4) and $\partial\tilde X=\emptyset$.
\endproclaim

\demo{Proof} We give the proof by induction on $n$. If $n=1$, then
$X$ is a circle. In this case, from (3.1.2), it is easy to see that
$f(X)$ is open in $\tilde X$ and that any point of $f(X)$ is an
interior one. Moreover, $f(X)$ is compact and closed in $\tilde X$
by Proposition 3.4. It then follows that $\tilde X=f(X)$, which is
compact and has an empty boundary (i.e. $\tilde X$ is also a
circle).

Now we assume that $n>1$. Since $\text{\rm
D}f:\Sigma_xX\to\Sigma_{\tilde x}\tilde X$ is also a
cone-neighborhood isometry for any $x\in X$ and $\tilde x\in f(x)$
(see Proposition 3.3), by the inductive assumption (note that
$\Sigma_pX$ is compact and has an empty boundary), we have that
$$\text{\rm D}f(\Sigma_xX)=\Sigma_{\tilde x}\tilde X\tag{3.7}$$
which is compact and has an empty boundary. It then remains to show
that $\tilde X=f(X)$. If this is not true, then, for an arbitrary
$\tilde q\in\tilde X\setminus f(X)$, there is a $\tilde p\in f(X)$
such that (note that $f(X)$ is compact by Proposition 3.4)
$$|\tilde q\tilde p|=\min\{|\tilde q\tilde p'|\ |\ \tilde p'\in f(X)\}.$$
Select a $[\tilde p\tilde q]$. By (1.1.2) and the first
variation formula, it is easy to see that there is an
$\varepsilon>0$ such that
$$V_{[\tilde p\tilde q],\varepsilon}\cap f(X)=\{\tilde p\},\tag{3.8}$$
where $V_{[\tilde p\tilde q],\varepsilon}$ is the $\varepsilon$-cone neighborhood of $[\tilde p\tilde q]$.
Let $p=f^{-1}(\tilde p)$. From the definition of $\text{\rm D}f:\Sigma_pX\to\Sigma_{\tilde p}\tilde
X$ (see Proposition 3.3 and its proof), (3.8) implies that
$$\uparrow_{\tilde p}^{\tilde q}\ \not\in\text{\rm D}f(\Sigma_pX),$$
which contradicts (3.7) (i.e., we obtain that $f(X)=\tilde X$).
\hfill$\qed$
\enddemo

From the proof of Theorem 3.5, we get the following corollary.

\proclaim{Corollary 3.6} Let $X, \tilde X\in \text{Alex}^n(1)$ with $n\geq1$, and let $f:X\rightarrow \tilde X$ be a
cone-neighborhood isometry. If $p\in X^\circ$, then $f(p)
\subset\tilde X^\circ$, and there is an $\epsilon>0$ such
that $B(\tilde p,\epsilon)\subseteq f(B(p,\epsilon))$ for any $\tilde p\in f(p)$.
\endproclaim

\remark{Remark \rm3.7} Let $f: X\to \tilde X$ be a cone-neighborhood
isometry. If $X$ is a Riemannian manifold, then $\tilde X=f(X)$ is a
Riemannian submanifold which may have more than one component. When
restricted to any component, $f^{-1}$  is a Riemannian covering map.
(Hint: since $\Sigma_xX$ is a unit sphere for any $x\in X$,
$\text{\rm D}f:\Sigma_xX\to\Sigma_{\tilde x}\tilde X$ is an
isometric embedding, which together with (3.1) implies that $x\in
X_m$.)

For general $X, \tilde X\in \text{Alex}(1)$ ($X$ is not a
Riemannian manifold), if $\dim(X)=\dim(\tilde X)$, then $f^{-1}$ is a
branch covering map (see Theorem 3.5). However, if $\dim(X)<\dim(\tilde X)$,
the above covering property of $f^{-1}$ may not hold.
\endremark

\vskip2mm

Nevertheless, the following property will be used in the proof of Theorem B.

\proclaim{Lemma 3.8} Let $f: X\to \tilde X$ ($X, \tilde X\in \text{Alex}(1)$) be a cone-neighborhood
isometry. If $X$ is compact and $\partial X=\emptyset$,
then
$$(f(X))^{\geq\frac\pi2}=(f(X))^{=\frac\pi2}. $$
\endproclaim

\demo{Proof} We give the proof by induction on $\dim(X)$.
If $\dim(X)=0$, then it is obvious
(because it is our convention that $X$ consists of
two points with distance $\pi$, so does
$f(X)$ by (3.1.2)). Now, we assume that $\dim(X)>0$. Let $\tilde y$
be an arbitrary point in $(f(X))^{\geq\frac\pi2}$, and let $\tilde
x_0\in f(X)$ such that $|\tilde y\tilde x_0|=\min\{|\tilde y\tilde
x||\ \tilde x\in f(X)\}$ (see Proposition 3.4). By the first variation
formula, for any $[\tilde x_0\tilde y]$ we have that
$\uparrow_{\tilde x_0}^{\tilde y}\in(\text{\rm
D}f(\Sigma_{x_0}X))^{\geq\frac\pi2}\subset\Sigma_{\tilde x_0}\tilde
X$ (see Proposition 3.3 for $\text{\rm D}f$), where
$x_0=f^{-1}(\tilde x_0)$. By induction, we have that $(\text{\rm
D}f(\Sigma_{x_0}X))^{\geq\frac\pi2}=(\text{\rm
D}f(\Sigma_{x_0}X))^{=\frac\pi2}$ in $\Sigma_{\tilde x_0}\tilde X$.
Hence, $[\tilde x_0\tilde y]$ is perpendicular to each segment of
$B(\tilde x_0,r_{\tilde x_0})\cap f([x_0z])$ (see (3.1)) for any
$[x_0z]\subset B(x_0,r_{\tilde x_0})$. By (1.1.2) and the choice of
$\tilde x_0$, it has to hold that $|\tilde y\tilde z|=\frac\pi2$ for
any $\tilde z\in B(\tilde x_0,r_{\tilde x_0})\cap f(X)$. It then is
not hard to conclude that $|\tilde y\tilde z|=\frac\pi2$ for any
$\tilde z\in f(X)$ (i.e. the lemma follows). \hfill$\qed$
\enddemo

We now begin to prove Key Lemma 1.7, starting with a special case.

\proclaim{Lemma 3.9} Key Lemma 1.7 is true in the case where $p_1\in X_1^m$
(see Key Lemma 1.6 for $X_1^m$).
\endproclaim

\demo{Proof}
We consider the multi-valued map
$$f_{p_1}:X_2\to(\Sigma_{p_1}X_1)^\perp \text{ by } p_2\mapsto\Uparrow_{p_1}^{p_2}.$$
By Key Lemma 1.6, Lemma 1.4 and (1.1.2), it is not hard to see that
$$f_{p_1}|_{X_2^\circ}:X_2^\circ\to(\Sigma_{p_1}X_1)^\perp \text{ is a cone-neighborhood isometry.}\tag{3.9}$$
Then the lemma follows from Corollary 3.6 and Theorem 3.5.
\hfill$\qed$
\enddemo

Based on Lemma 3.9, we can conclude another important property of $X_i^m$.

\proclaim{Lemma 3.10} For any $p_2\in X_2^\circ$, there exists
a positive number $m_{p_2}$ ($\leq m$) such that $\lambda_{p_1p_2}=m_{p_2}$ \
for all $p_1\in X_1^m$ and $\lambda_{p_1'p_2}\leq m_{p_2}$ for all
$p_1'\in X_1^\circ$. (And the similar statement holds for $p_1\in X_1^\circ$.)
\endproclaim

\demo{Proof} By Lemma 3.9, for $p_1\in X_1^m$,  any
$\uparrow_{p_1}^{p_2}$ belongs to $((\Sigma_{p_1}X_1)^\perp)^\circ$.
Then (2.11) together with Lemma 1.4 implies that
$$\lambda_{p_2p_1}\geq \lambda_{p_2p_1'} \text{ for all } p_1'\in X_1^\circ,$$
which implies the lemma.
\hfill$\qed$
\enddemo

Now, we can complete the proof of Key Lemma 1.7.

\demo{Proof of Key Lemma 1.7}

The proof is almost identical to that
of Lemma 3.9 except that we should use Lemma 3.10 instead of Key
Lemma 1.6 to conclude that
$f_{p_1}|_{X_2^\circ}:X_2^\circ\to(\Sigma_{p_1}X_1)^\perp$ is a
cone-neighborhood isometry. \hfill$\qed$
\enddemo

\vskip4mm

%%%%%%%%%%%%%%%%%%%%%%%%%%%%%%%%%%% Section 4 %%%%%%%%%%%%%%%%%%%%%%%%%%%%%%%%%%%%%%%%%%%

\head 4. Proof of Theorem A
\endhead

\vskip4mm

To prove Theorem A, by Lemma 1.2 and Theorem 1.8 it remains to prove
(A2), where $X_i$ is compact and has an empty boundary. As outlined
at the end of Introduction, for any $p_i\in X_i^m$ (see Key Lemma
1.6), we shall construct a finite group $\Gamma$ ($|\Gamma|=m$)
acting isometrically on $(\Sigma_{p_i}X_i)^\perp$ such that
$(\Sigma_{p_i}X_i)^\perp/\Gamma \overset{\text{isom}}\to\cong X_j$
($i\neq j\in \{1,2\}$); and based on Key Lemma 1.7 we can define a
natural map, $((\Sigma_{p_1}X_1)^\perp*(\Sigma_{p_2}X_2)^\perp
)/\Gamma\to X$, and check that the map is an isometry.

According to Key Lemma 1.7 and its proof, for any $p_i\in X_i$, we
have that
$$f_{p_i}:X_j\to(\Sigma_{p_i}X_i)^\perp \text{ by } p_j\mapsto\Uparrow_{p_i}^{p_j}\tag{4.1}$$
is a cone-neighborhood isometry ($i\neq j$), and that
$\Uparrow_{p_i}^{X_j}=(\Sigma_{p_i}X_i)^\perp$ which is compact and
has an empty boundary.  It follows that  $$\text{\it For any
$\zeta\in (\Sigma_{p_i}X_i)^\perp$, there is a $[p_ip_j]$ with
$p_j\in X_j$ such that $\zeta=\uparrow_{p_i}^{p_j}$}.\tag{4.2}$$

Using the above consequences of Key Lemma 1.7, we can explore a local
join structure on $X$.

\proclaim{Lemma 4.1} {\rm (4.1.1)} Any $x\in X$ belongs to some
$[p_1p_2]$ with $p_i\in X_i$.

\vskip1mm

\noindent{\rm(4.1.2)} For any $p_i\in X_i^m$,
$\Sigma_{p_i}X=\Sigma_{p_i}X_i*(\Sigma_{p_i}X_i)^\perp$.

\endproclaim

\demo{Proof} (4.1.1) Since $X_1$ is compact, there is a $p_1\in X_1$
such that $|xp_1|=\min\{|xp_1'||\ p_1'\in X_1\}$. By Lemma 1.3, we
have that $|xp_1|\leq\frac\pi2$; and by the first variation formula,
for any $[p_1x]$, we know that $\uparrow_{p_1}^x$ belongs to
$(\Sigma_{p_1}X_1)^{\geq\frac\pi2}$ in $\Sigma_{p_1}X$. Again by
Lemma 1.3, $\uparrow_{p_1}^x\in (\Sigma_{p_1}X_1)^\perp$
($=(\Sigma_{p_1}X_1)^{=\frac\pi2}$), so by (4.2) there exists
a $[p_1p_2]$ with $p_2\in X_2$ such that
$\uparrow_{p_1}^x=\uparrow_{p_1}^{p_2}$. This together with
``$|xp_1|\leq\frac\pi2$'' implies that $[p_1x]$ ($\ni x$)
lies in $[p_1p_2]$.

\vskip1mm

(4.1.2) Since $(\Sigma_{p_i}X_i)^\perp$ has an empty boundary  and
$p_i\in X_i^m$, we conclude that
$\Sigma_{p_i}X=\Sigma_{p_i}X_i*(\Sigma_{p_i}X_i)^\perp$ by Corollary
2.5. \hfill$\qed$
\enddemo

In view of the above, the join structure on $X$ should be independent of
$p_i\in X_i^m$, which will be justified in Lemma 4.2.

Let $p_1, p_1'\in X_1^m$, and $[p_1p_1']\subset X_1$. Since $\Sigma_{p_1}X=\Sigma_{p_1}X_1*(\Sigma_{p_1}X_1)^\perp$ (resp.
$\Sigma_{p_1'}X=\Sigma_{p_1'}X_1*(\Sigma_{p_1'}X_1)^\perp$) by (4.1.2),
for any $p_2\in X_2$ and given $[p_1p_2]$ (resp. $[p_1'p_2]$),
by Lemma 1.4 there is a unique $[p_1'p_2]$ (resp. $[p_1p_2]$)
such that $[p_1p_2]$, $[p_1'p_2]$ and $[p_1p_1']$ bound
a convex spherical surface. Hence, due to
(4.2), we can define a 1-1 map
$$h_{[p_1p_1']}:(\Sigma_{p_1}X_1)^\perp\to (\Sigma_{p_1'}X_1)^\perp
\ \text{ by }\uparrow_{p_1}^{p_2}\mapsto\uparrow_{p_1'}^{p_2}.$$

\proclaim{Lemma 4.2} The map $h_{[p_1p_1']}:(\Sigma_{p_1}X_1)^\perp
\to(\Sigma_{p_1'}X_1)^\perp$ is an isometry.
\endproclaim

\demo{Proof} For convenience, we let $h$ denote $h_{[p_1p_1']}$ in this proof.
We need to show that
$$|h(\uparrow_{p_1}^{p_2})h(\uparrow_{p_1}^{p_2'})|=|\uparrow_{p_1}^{p_2}\uparrow_{p_1}^{p_2'}|\tag{4.3}$$
for any $\uparrow_{p_1}^{p_2}, \uparrow_{p_1}^{p_2'}\in
(\Sigma_{p_1}X_1)^\perp$. Let $\tilde\gamma(t)|_{t\in[0,1]}$ be a
minimal geodesic in $(\Sigma_{p_1}X_1)^\perp$ with
$\tilde\gamma(0)=\uparrow_{p_1}^{p_2}$ and
$\tilde\gamma(1)=\uparrow_{p_1}^{p_2'}$. By (4.1) and (4.2), we can
consider $\gamma(t)\triangleq f_{p_1}^{-1}(\tilde\gamma(t))$ for all
$t\in [0,1]$ (note that $\gamma(0)=p_2$ and $\gamma(1)=p_2'$).

\vskip1mm

\noindent{\bf Claim 1}: {\it $\gamma(t)|_{t\in[0,1]}$ is a piecewise geodesic in $X_2$
with length equal to that of $\tilde \gamma(t)|_{t\in[0,1]}$, $|\uparrow_{p_1}^{p_2}\uparrow_{p_1}^{p_2'}|$.}
This is due to that $f_{p_1}:X_2\rightarrow (\Sigma_{p_1}X_1)^\perp$
is a cone-neighbor-hood
isometry and that $f_{p_1}$ is surjective (see Remark 3.2 and (4.1)).

\vskip1mm

\noindent{\bf Claim 2}: {\it $h$ is a continuous map.}
We need to show that if $\uparrow_{p_1}^{p_2'}\to \uparrow_{p_1}^{p_2}$
(in $(\Sigma_{p_1}X_1)^\perp$), then $h(\uparrow_{p_1}^{p_2'})\to h(\uparrow_{p_1}^{p_2})$.
Let $\uparrow_{p_1}^{p_2}$, $\uparrow_{p_1}^{p_2'}$, $h(\uparrow_{p_1}^{p_2})$
and $h(\uparrow_{p_1}^{p_2'})$ be the directions of
$[p_1p_2]$, $[p_1p_2']$, $[p_1'p_2]$ and $[p_1'p_2']$ respectively (see (4.2)).
Note that $\Sigma_{p_1}X=\Sigma_{p_1}X_1*(\Sigma_{p_1}X_1)^\perp$,
so the (unique) $[\uparrow_{p_1}^{p_2'}\uparrow_{p_1}^{p_1'}]$
converges to the $[\uparrow_{p_1}^{p_2}\uparrow_{p_1}^{p_1'}]$ as
$\uparrow_{p_1}^{p_2'}\to \uparrow_{p_1}^{p_2}$.
This implies that the convex
spherical surface $\overline{\triangle p_2'p_1p_1'}$ (which realizes $[\uparrow_{p_1}^{p_2'}\uparrow_{p_1}^{p_1'}]$)
bounded by $[p_1p_2']$, $[p_1'p_2']$ and $[p_1p_1']$
converges to
$\overline{\triangle p_2p_1p_1'}$ (which realizes $[\uparrow_{p_1}^{p_2}\uparrow_{p_1}^{p_1'}]$) bounded by $[p_1p_2]$, $[p_1'p_2]$ and $[p_1p_1']$,
so $h(\uparrow_{p_1}^{p_2'})$ converges to $h(\uparrow_{p_1}^{p_2})$.

\vskip1mm

By Claim 2, $h(\tilde\gamma(t))|_{t\in[0,1]}\subset f_{p_1'}(\gamma(t))|_{t\in[0,1]}$
is a continuous curve in $(\Sigma_{p_1'}X_1)^\perp$, while
$\gamma(t)|_{t\in[0,1]}$ is a piecewise geodesic in $X_2$ by Claim 1.
Note that $f_{p_1'}:X_2\rightarrow (\Sigma_{p_1'}X_1)^\perp$
(defined by $p_2''\mapsto\Uparrow_{p_1'}^{p_2''}$) is also a cone-neighborhood
isometry. Hence, $h(\tilde\gamma(t))|_{t\in[0,1]}$
is also a piecewise geodesic with length equal to that of $\gamma(t)|_{t\in[0,1]}$.
It then follows that
$$|h(\uparrow_{p_1}^{p_2})h(\uparrow_{p_1}^{p_2'})|\leq|\uparrow_{p_1}^{p_2}\uparrow_{p_1}^{p_2'}|.$$
In turn, we can similarly get that
$$|h(\uparrow_{p_1}^{p_2})h(\uparrow_{p_1}^{p_2'})|\geq
|h^{-1}(h(\uparrow_{p_1}^{p_2}))h^{-1}(h(\uparrow_{p_1}^{p_2'}))|=|\uparrow_{p_1}^{p_2}\uparrow_{p_1}^{p_2'}|.$$
Hence, (4.3) follows.
\hfill$\qed$
\enddemo

\remark{Remark \rm4.3} From the proof of Lemma 4.2, it is not hard
to see that there are neighborhoods $U_{p_2}$,
$U_{\uparrow_{p_1}^{p_2}}$ ($\subset (\Sigma_{p_1}X_1)^\perp$) and
$U_{\uparrow_{p_1'}^{p_2}}$ ($\subset (\Sigma_{p_1'}X_1)^\perp$) of
$p_2$, $\uparrow_{p_1}^{p_2}$ and $\uparrow_{p_1'}^{p_2}$
respectively such that the following diagram commutes
$$\CD
 U_{p_2} & @>f_{p_1}>> & U_{\uparrow_{p_1}^{p_2}}\\
   @V\text{id}VV &                 & @VV h_{[p_1p_1']}\ , V\\
 U_{p_2}    & @>f_{p_1'}>> & U_{\uparrow_{p_1'}^{p_2}}
\endCD\hskip0.7cm $$
\noindent in which each map is an isometry.
\endremark

Based on Lemma 4.2, we will prove the following property
that is crucial to our construction of $\Gamma$.

\proclaim{Lemma 4.4} For $p_i\in X_i^m$, any two minimal geodesics
$[p_1p_2]_1$ and $[p_1p_2]_2$ determine an isometry, $\hat\gamma:
(\Sigma_{p_1}X_1)^\perp\to (\Sigma_{p_1}X_1)^\perp$, such that $\hat\gamma((\uparrow_{p_1}^{p_2})_1)=(\uparrow_{p_1}^{p_2})_2$ and
$\hat\gamma(\Uparrow_{p_1}^{p_2'})=\Uparrow_{p_1}^{p_2'}$ for
$p_2'\in X_2$.
\endproclaim

Before giving a proof, we need some preparations.

\proclaim{Lemma 4.5} For $p_1, p_1'\in X_1$, $p_2\in X_2$ and given
$[p_1p_2]$ and $[p_1p_1']$, any $[\uparrow_{p_1}^{p_2}\uparrow_{p_1}^{p_1'}]$ in
$\Sigma_{p_1}X$ is realized by a convex spherical surface bounded by
a triangle $\triangle p_2p_1p_1'$ containing $[p_1p_2]$ and
$[p_1p_1']$. (And the similar statement holds for given
$[p_2p_1]$ and $[p_2p_2']\subset X_2$.)
\endproclaim

\demo{Proof} Select an arbitrary $[\uparrow_{p_1}^{p_2}\uparrow_{p_1}^{p_1'}]$,
and let $\eta\in[\uparrow_{p_1}^{p_2}\uparrow_{p_1}^{p_1'}]^\circ$. We first
show that the lemma is true if there is a $[p_1q]$ such
that $\uparrow_{p_1}^q=\eta$.

{\bf Claim}: {\it $[p_1q]$ lies in a convex spherical
surface bounded by a triangle $\triangle p_1\bar p_1\bar p_2$ with
$\bar p_i\in X_i$}. By (4.1.1), there is a $[\bar p_1\bar p_2]$ with
$\bar p_i\in X_i$ such that $q\in [\bar p_1\bar p_2]$. By Lemma 1.4,
there is a convex spherical surface containing $p_1$ and $[\bar
p_1\bar p_2]$. Then by (1.1.3), $[p_1q]$ and $[\bar p_1\bar p_2]$
determine a triangle $\triangle p_1\bar p_1\bar p_2$ bounding a
convex spherical surface $D$ (containing $[p_1q]$ and $[\bar p_1\bar
p_2]$).

We assume that the $\triangle p_1\bar p_1\bar p_2$ in the claim is
formed by $[p_1\bar p_2]$, $[p_1\bar p_1]$ and $[\bar p_1\bar p_2]$.
It then follows from the claim that
$\uparrow_{p_1}^q\in[\uparrow_{p_1}^{p_2}\uparrow_{p_1}^{p_1'}]\cap[\uparrow_{p_1}^{\bar
p_2}\uparrow_{p_1}^{\bar p_1}]$, where $[\uparrow_{p_1}^{\bar
p_2}\uparrow_{p_1}^{\bar p_1}]\ (\subset\Sigma_{p_1}X)$ is realized
by $D$. Note that
$|\uparrow_{p_1}^{p_2}\uparrow_{p_1}^{p_1'}|=|\uparrow_{p_1}^{p_2}\uparrow_{p_1}^{\bar
p_1}| =|\uparrow_{p_1}^{\bar p_2}\uparrow_{p_1}^{\bar
p_1}|=|\uparrow_{p_1}^{\bar p_2}\uparrow_{p_1}^{p_1'}|=\frac\pi2$ in
$\Sigma_{p_1}X$, so it has to hold that $[\uparrow_{p_1}^{\bar
p_2}\uparrow_{p_1}^{\bar
p_1}]=[\uparrow_{p_1}^{p_2}\uparrow_{p_1}^{p_1'}]$, which implies
that $[p_1p_2]=[p_1\bar p_2]$ and $[p_1\bar p_1]\supseteq[p_1p_1']$
or $[p_1\bar p_1]\subset[p_1p_1']$. If $[p_1\bar
p_1]\supseteq[p_1p_1']$, then the proof is done. If $[p_1\bar
p_1]\subset[p_1p_1']$, then, still by Lemma 1.4, there is a
$[p_2p_1']$ which together with $[p_2\bar p_1]\ (=[\bar p_2\bar
p_1])$ and $[\bar p_1p_1']$ bounds a convex spherical surface $D'$.
It is not hard to see that $D\cup D'$ is a convex spherical surface
(bounded by $[p_1p_2]$, $[p_1p_1']$ and $[p_1'p_2]$) which realizes
$[\uparrow_{p_1}^{p_2}\uparrow_{p_1}^{p_1'}]$.

Next, it suffices to show that there must be a
$[p_1q]$ such that $\uparrow_{p_1}^q=\eta$ for any
$\eta\in[\uparrow_{p_1}^{p_2}\uparrow_{p_1}^{p_1'}]^\circ$. If there is no $[p_1q]$ such that
$\uparrow_{p_1}^q=\eta$, we can select
$\{[p_1q_j]\}_{j=1}^\infty$ with $|p_1q_j|\to0$ as $j\to\infty$ such
that $\uparrow_{p_1}^{q_j}\to\eta$ ([BGP]). By the
claim above, $[p_1q_j]$ lies in a convex spherical surface $D_j$ bounded by
a triangle $\triangle p_1p_1^jp_2^j$ with $p_i^j\in X_i$.
Assume that the $\triangle p_1p_1^jp_2^j$ is formed
by $[p_1p_2^j]$, $[p_1p_1^j]$ and $[p_1^jp_2^j]$.
Note that $D_j$ realizes a $[\uparrow_{p_1}^{p_2^j}\uparrow_{p_1}^{p_1^j}]$ in $\Sigma_{p_1}X$ with
$\uparrow_{p_1}^{q_j}\in
[\uparrow_{p_1}^{p_2^j}\uparrow_{p_1}^{p_1^j}]$, so it follows from
$\uparrow_{p_1}^{q_j}\to\eta$ that
$$[\uparrow_{p_1}^{p_2^j}\uparrow_{p_1}^{p_1^j}]\to
[\uparrow_{p_1}^{p_2}\uparrow_{p_1}^{p_1'}]\tag{4.4}$$
(similar to ``$[\uparrow_{p_1}^{\bar p_2}\uparrow_{p_1}^{\bar
p_1}]=[\uparrow_{p_1}^{p_2}\uparrow_{p_1}^{p_1'}]$'' in the above paragraph),
which implies that $[p_1p_2^j]\to [p_1p_2]$. Now we consider the cone-neighborhood isometry $f_{p_2^j}:X_1\to
(\Sigma_{p_2^j}X_2)^\perp$ (see (4.1)). Note that $D_j$ realizes a
$[\uparrow_{p_2^j}^{p_1}\uparrow_{p_2^j}^{p_1^j}]$ which belongs to
$f_{p_2^j}([p_1p_1^j])$ in $(\Sigma_{p_2^j}X_2)^\perp$. By Claim 2 in the
proof of Proposition 3.3, there is a convex spherical surface $\bar D_j$
bounded by a triangle $\triangle p_2^jp_1p_1'$ containing the
$[p_1p_2^j]$ and $[p_1p_1']$ such that $\bar D_j$ realizes a
$[\uparrow_{p_2^j}^{p_1}\uparrow_{p_2^j}^{p_1'}]$
(belonging to $f_{p_2^j}([p_1p_1'])$ in $(\Sigma_{p_2^j}X_2)^\perp$) with
$$\left|\xuparrow{0.5cm}_{\uparrow_{p_2^j}^{p_1}}^{\uparrow_{p_2^j}^{p_1^j}}
\xuparrow{0.5cm}_{\uparrow_{p_2^j}^{p_1}}^{\uparrow_{p_2^j}^{p_1'}}\right|=
\left|\uparrow_{p_1}^{p_1^j}\uparrow_{p_1}^{p_1'}\right|.$$
This together with
``$\uparrow_{p_1}^{p_1^j}\to \uparrow_{p_1}^{p_1'}$'' (see (4.4)) implies that
the $[\uparrow_{p_1}^{p_2^j}\uparrow_{p_1}^{p_1'}]$ in
$\Sigma_{p_1}X$ realized by $\bar D_j$ is sufficiently close to
$[\uparrow_{p_1}^{p_2^j}\uparrow_{p_1}^{p_1^j}]$. It then follows from (4.4) that
$[\uparrow_{p_1}^{p_2^j}\uparrow_{p_1}^{p_1'}]\to
[\uparrow_{p_1}^{p_2}\uparrow_{p_1}^{p_1'}]$, so
$\bar D_j$ converges to a convex spherical surface $D$ bounded by a
triangle $\triangle p_2p_1p_1'$ containing the $[p_1p_2]$ and
$[p_1p_1']$ such that the
$[\uparrow_{p_1}^{p_2}\uparrow_{p_1}^{p_1'}]$ is realized by $D$.
This contradicts the assumption that there is no $[p_1q]$
such that $\uparrow_{p_1}^q=\eta$. \hfill$\qed$
\enddemo

Based on Lemma 4.5, we have a further observation.

\proclaim{Lemma 4.6} Let $p_i\in X_i^m$ ($i=1,2$). For any $p_1'\in
X_1$, $m_{p_1'}$ (see Lemma 3.10 for $m_{p_1'}$) divides $m$; and
given $[p_1'p_2]$ and $[p_1'p_1]$, there are $\frac{m}{m_{p_1'}}$ many
minimal geodesics in $\Sigma_{p_1'}X$ from $\uparrow_{p_1'}^{p_2}$
to $\uparrow_{p_1'}^{p_1}$.
\endproclaim

\demo{Proof} At first, note that $\lambda_{p_2p_1'}=m_{p_1'}$ by
Lemma 3.10. Next, note that (1.5) and Remark 1.9 enable us to
apply Key Lemma 1.6 and Lemma 3.10 on $\Sigma_{p_1'}X$
to conclude that there is an open and dense subset $A$ in
$(\Sigma_{p_1'}X_1)^\perp$ such that for any fixed
$\xi\in\Sigma_{p_1'}X_1$ there is an $m_\xi$ such that
$$\lambda_{\xi\zeta}=m_\xi\ \forall\ \zeta\in A \text{ and }
\lambda_{\xi\zeta'}\leq m_\xi\ \forall\ \zeta'\in
(\Sigma_{p_1'}X_1)^\perp.$$ At last, note that, for any fixed
$[p_1p_2]$, there is a unique $[p_1'p_2]$ which together with
$[p_1p_2]$ and $[p_1p_1']$ bounds a convex spherical surface (by
Lemma 1.4 and (4.1.2)). Then by Lemma 4.5, it suffices to show
that, for any fixed $[p_1'p_2]$, there are
$m_{\uparrow_{p_1'}^{p_1}}$ pieces of convex spherical surfaces
bounded by triangles containing  $[p_1'p_1]$ and $[p_1'p_2]$ (which
implies that $m=m_{p_1'}\cdot m_{\uparrow_{p_1'}^{p_1}}$).

Note that we can find $\zeta_j\in A$ such that $\zeta_j\to\uparrow_{p_1'}^{p_2}$
as $j\to\infty$, and select $[p_1'p_2^j]$ with $p_2^j\in X_2$
such that $\uparrow_{p_1'}^{p_2^j}=\zeta_j$ (see (4.2)). By Lemma 4.5,
there are $[p_1p_2^j]_k$ with $1\leq k\leq m_{\uparrow_{p_1'}^{p_1}}$ such that
the minimal geodesics between $\uparrow_{p_1'}^{p_1}$ and $\zeta_j$ are realized
by convex spherical surfaces $D_k^j$ bounded
by $[p_1p_1']$, $[p_1'p_2^j]$ and $[p_1p_2^j]_k$.
Since $p_2\in X_2^m$, we can assume that $p_2^j\in X_2^m$ ($X_2^m$ is open).
Then by (2.4) (and Lemma 1.4), each $[p_1p_2^j]_k$ converges to some $[p_1p_2]_k$ as $j\to \infty$
with $[p_1p_2]_{k_1}\neq[p_1p_2]_{k_2}$ for $1\leq k_1\neq k_2\leq m_{\uparrow_{p_1'}^{p_1}}$.
This implies that
each $D_k^j$ converges to a convex spherical surface
bounded by $[p_1'p_1]$, $[p_1'p_2]$ and $[p_1p_2]_k$, so the proof is done.
\hfill$\qed$
\enddemo

We are ready to give a proof of Lemma 4.4.

\demo{Proof of Lemma 4.4}

Let $\tilde\gamma(t)|_{t\in[0,1]}$ be a
minimal geodesic in $(\Sigma_{p_2}X_2)^\perp$ with
$\tilde\gamma(0)=(\uparrow_{p_2}^{p_1})_1$ and
$\tilde\gamma(1)=(\uparrow_{p_2}^{p_1})_2$; and let
$\gamma(t)|_{[0,1]}\triangleq
f_{p_2}^{-1}(\tilde\gamma(t)|_{[0,1]})$ (note that
$\gamma(0)=\gamma(1)=p_1$),  which is a piecewise geodesic in $X_1$
with length equal to
$|(\uparrow_{p_2}^{p_1})_1(\uparrow_{p_2}^{p_1})_2|$ (see Claim 1 in the proof
of Lemma 4.2).

Case 1: $\gamma(t)|_{[0,1]}\subset X_1^m$.  Since
$\gamma(t)|_{[0,1]}$ is a piecewise geodesic, we select
$0=t_0<t_1<\cdots<t_l=1$ such that $\gamma(t)|_{[t_k,t_{k+1}]}$ is a
minimal geodesic for $k=0,\cdots, l-1$. By Lemma 4.2, we know that
$\gamma(t)|_{[t_k,t_{k+1}]}$ determines an isometry
$g_{\gamma(t)|_{[t_k,t_{k+1}]}}:(\Sigma_{\gamma(t_k)}X_1)^\perp\to
(\Sigma_{\gamma(t_{k+1})}X_1)^\perp.$ Composing them, we get an
isometry
$$\hat\gamma\triangleq g_{\gamma(t)|_{[t_{l-1},t_l]}}\circ\cdots\circ
g_{\gamma(t)|_{[t_1,t_2]}}\circ g_{\gamma(t)|_{[t_0t_{1}]}}:
(\Sigma_{p_1}X_1)^\perp\to (\Sigma_{p_1}X_1)^\perp.\tag{4.5}$$ It is
easy to see that
$\hat\gamma((\uparrow_{p_1}^{p_2})_1)=(\uparrow_{p_1}^{p_2})_2$, and
$\hat\gamma(\uparrow_{p_1}^{p_2'})\in\Uparrow_{p_1}^{p_2'}$ for any
$\uparrow_{p_1}^{p_2'}$ with $p_2'\in X_2$. Moreover, due to Remark
4.3, it is not hard to see that there are neighborhoods $U_{p_2}$,
$U_{(\uparrow_{p_1}^{p_2})_1}$ and $U_{(\uparrow_{p_1}^{p_2})_2}$ of
$p_2$, $(\uparrow_{p_1}^{p_2})_1$ and $(\uparrow_{p_1}^{p_2})_2$
respectively such that the following diagram commutes
$$\CD
 U_{p_2} & @>f_{p_1}>> & U_{(\uparrow_{p_1}^{p_2})_1}\\
   @V\text{id}VV &                 & @VV \hat\gamma V\\
 U_{p_2}    & @>f_{p_1}>> & U_{(\uparrow_{p_1}^{p_2})_2}
\endCD.\tag{4.6} $$
Hence, $\hat\gamma$ does not depend on the choice of
$\{t_k\}_{k=0}^l$. (In fact, for all piecewise geodesics
$\tilde\gamma(t)|_{t\in[0,1]}\subset(\Sigma_{p_2}{X_2})^\perp$
between $(\uparrow_{p_2}^{p_1})_1$ and $(\uparrow_{p_2}^{p_1})_2$
with $f_{p_2}^{-1}(\tilde\gamma(t))\in X_1^m$, the isometries of
$(\Sigma_{p_1}{X_1})^\perp$ gotten through such an way are the same
one.)

\vskip1mm

Case 2: There is $t_0\in(0,1)$ such that $\gamma(t_0)\not\in X_1^m$,
and there is no other $[p_2p_1]$ such that
$|\uparrow_{p_2}^{p_1}\tilde\gamma(t_0)|\leq|(\uparrow_{p_2}^{p_1})_j\tilde\gamma(t_0)|$
for $j=1$ and 2. In this case, by Lemma 4.6
(and its proof), it is not hard to see that
both $\gamma(t)|_{[0,t_0]}$ and $\gamma(t)|_{[t_0,1]}$ are a same
minimal geodesic between $p_1$ and $\gamma(t_0)$, denoted by
$[p_1\gamma(t_0)]$, and
that there is an open and dense subset $A$ in
$(\Sigma_{\gamma(t_0)}X_1)^\perp$ such that
$$\lambda_{\uparrow_{\gamma(t_0)}^{p_1}\zeta}=2\ \forall\ \zeta\in A \text{ and }
\lambda_{\uparrow_{\gamma(t_0)}^{p_1}\zeta'}\leq 2\ \forall\
\zeta'\in (\Sigma_{\gamma(t_0)}X_1)^\perp. \tag{4.7}$$ By Lemma 1.4
and (4.1.2), for any fixed $[p_1p_2']$ with $p_2'\in X_2$, there is
a unique $[\gamma(t_0)p_2']$ which together with $[p_1p_2']$ and
$[p_1\gamma(t_0)]$ bounds a convex spherical surface $D$. Then (4.7)
implies that, besides $[p_1p_2']$, there exists at most one $[p_1p_2']'$
which together with $[\gamma(t_0)p_2']$ and
$[p_1\gamma(t_0)]$ bounds a convex spherical surface $D'$. Then we
can define a 1-1 map
$$\hat\gamma:(\Sigma_{p_1}X_1)^\perp\to (\Sigma_{p_1}X_1)^\perp \text{ by }
\uparrow_{p_1}^{p_2'}\mapsto (\uparrow_{p_1}^{p_2'})' \text{ and }
(\uparrow_{p_1}^{p_2'})'\mapsto \uparrow_{p_1}^{p_2'},\tag{4.8}$$
where $(\uparrow_{p_1}^{p_2'})'=\uparrow_{p_1}^{p_2'}$ if
$[p_1p_2']'$ does not exist. Note that
$\hat\gamma((\uparrow_{p_1}^{p_2})_1)=(\uparrow_{p_1}^{p_2})_2$ and
$\hat\gamma(\uparrow_{p_1}^{p_2'})$ $\in\Uparrow_{p_1}^{p_2'}$. Similar to
the proof of Lemma 4.2, we will conclude that $\hat\gamma$ is an
isometry once we show that it is a continuous map. Hence, it
suffices to show that $[p_1p_2^j]'\to[p_1p_2']'$ if
$[p_1p_2^j]\to[p_1p_2']$ as $j\to\infty$ (where $p_2^j\in X_2$). Let
$[\gamma(t_0)p_2^j]$ be the unique minimal geodesic which together
with $[p_1p_2^j]$ (resp. $[p_1p_2^j]'$) and $[p_1\gamma(t_0)]$
bounds a convex spherical surface $D_j$ (resp. $D_j'$). Note that
$D_j$ (resp. $D_j'$) realizes a minimal geodesic
$[\uparrow_{\gamma(t_0)}^{p_2^j}\uparrow_{\gamma(t_0)}^{p_1}]$
(resp.
$[\uparrow_{\gamma(t_0)}^{p_2^j}\uparrow_{\gamma(t_0)}^{p_1}]'$) in
$\Sigma_{\gamma(t_0)}X$. Due to (4.1.2), we have that
$D_j\to D $ as $j\to\infty$ because $[p_1p_2^j]\to[p_1p_2']$, so
$[\uparrow_{\gamma(t_0)}^{p_2^j}\uparrow_{\gamma(t_0)}^{p_1}]$
converges to
$[\uparrow_{\gamma(t_0)}^{p_2'}\uparrow_{\gamma(t_0)}^{p_1}]$ which
is realized by $D$. Then according to (4.7) and Lemma 1.4 on $\Sigma_{\gamma(t_0)}X$,
$[\uparrow_{\gamma(t_0)}^{p_2^j}\uparrow_{\gamma(t_0)}^{p_1}]'$ has
to converge to
$[\uparrow_{\gamma(t_0)}^{p_2'}\uparrow_{\gamma(t_0)}^{p_1}]'$ which
is realized by $D'$. It follows that $D_j'\to D'$ as $j\to\infty$,
so $[p_1p_2^j]'\to[p_1p_2']'$. (Similarly, $\hat\gamma$
satisfies a corresponding diagram like (4.6).)

\vskip1mm

Case 3: There is $t_0\in(0,1)$ such that $\gamma(t_0)\not\in X_1^m$,
and there is a $[p_2p_1]_3$ such that
$|(\uparrow_{p_2}^{p_1})_3\tilde\gamma(t_0)|\leq|(\uparrow_{p_2}^{p_1})_j\tilde\gamma(t_0)|$
for $j=1$ and 2. Since $\tilde\gamma(t)|_{t\in[0,1]}$ is
a minimal geodesic in $(\Sigma_{p_2}X_2)^\perp$ between $(\uparrow_{p_2}^{p_1})_1$ and
$(\uparrow_{p_2}^{p_1})_2$, it has to hold that
$$|(\uparrow_{p_2}^{p_1})_1(\uparrow_{p_2}^{p_1})_3|<|(\uparrow_{p_2}^{p_1})_1(\uparrow_{p_2}^{p_1})_2|
\text{ and }
|(\uparrow_{p_2}^{p_1})_2(\uparrow_{p_2}^{p_1})_3|<|(\uparrow_{p_2}^{p_1})_1(\uparrow_{p_2}^{p_1})_2|.$$
Then we can first look for the isometry of $(\Sigma_{p_1}X_1)^\perp$ determined by $[p_1p_2]_1$
and $[p_1p_2]_3$ (and $[p_1p_2]_2$ and $[p_1p_2]_3$) by repeating the above process.
Since there are only $m$ pieces of minimal geodesics between
$p_1$ and $p_2$, we can find the wanted isometry $\hat\gamma$ of
$(\Sigma_{p_1}X_1)^\perp$, a composition of some isometries in Case 1 or 2.
Moreover, the $\hat\gamma$ satisfies a corresponding diagram like (4.6).
\hfill$\qed$
\enddemo

In the rest of this section, we will complete the proof of (A2) in Theorem A.

\demo{Proof of (A2) in Theorem A}

If $\dim(X_1)=0$ (or $\dim(X_2)=0$), then, by our convention,
$\partial X_1=\emptyset$ means that $X_1$ consists of two points
with distance $\pi$. This implies that $X=X_1*X_2$. From now on, we
will assume that $\dim(X_i)>0$ for $i=1$ and 2.

\vskip1mm

We first fix two arbitrary points $p_1\in X_1^m$ and $p_2\in X_2^m$,
and let $\{[p_1p_2]_k\}_{k=1}^m$ be all minimal geodesics between
$p_1$ and $p_2$. According to Lemma 4.4, $[p_1p_2]_j$ and
$[p_1p_2]_h$ determines an isometry of $(\Sigma_{p_1}X_1)^\perp$,
denoted by $\gamma_{jh}$, with
$\gamma_{jh}((\uparrow_{p_1}^{p_2})_j)=(\uparrow_{p_1}^{p_2})_h$.
And from the proof of Lemma 4.4 (especially by (4.6)), we know that
$$\gamma_{hl}\cdot\gamma_{jh}=\gamma_{jl};$$
and because
$\gamma_{jh}((\uparrow_{p_1}^{p_2})_1)=(\uparrow_{p_1}^{p_2})_k$ for
some $k$, we have that $\gamma_{jh}=\gamma_{1k}$. Hence,
$$\Gamma_2\triangleq\{\gamma_{11},
\gamma_{12},\cdots,\gamma_{1m}\}$$ is a finite group (where
$\gamma_{11}$ is the unit element) which acts on
$(\Sigma_{p_1}X_1)^\perp$ by isometries; and the orbit of any
$\uparrow_{p_1}^{p_2'}$ (with $p_2'\in X_2$) under this action is
just $\Uparrow_{p_1}^{p_2'}$. Based on Lemma 1.4, it is easy to see
that the following map
$$\pi_2:(\Sigma_{p_1}X_1)^\perp/\Gamma_2\longrightarrow X_2 \text{ by }
\Uparrow_{p_1}^{p_2'}\longmapsto p_2'$$ is an isometry.

Similarly, $\{[p_1p_2]_k\}_{k=1}^m$ also determines a finite group
$$\Gamma_1\triangleq\{\bar\gamma_{11},
\bar\gamma_{12},\cdots,\bar\gamma_{1m}\}$$ which acts on
$(\Sigma_{p_2}X_2)^\perp$ by isometries with
$\bar\gamma_{jh}((\uparrow_{p_2}^{p_1})_j)=(\uparrow_{p_2}^{p_1})_h$;
and $$\pi_1:(\Sigma_{p_2}X_2)^\perp/\Gamma_1\longrightarrow X_1
\text{ by } \Uparrow_{p_2}^{p_1'}\longmapsto p_1'\tag{4.9}$$ is also
an isometry.

{\bf Claim 1}: {\it The map $$g:\Gamma_2\to \Gamma_1 \text{ defined by
} \gamma_{1k}\mapsto \bar\gamma_{1k}^{-1}$$ is an isomorphism}. We
need only to show that $g$ is a homomorphism, i.e.,
$$g(\gamma_{1k}\cdot\gamma_{1k'})=g(\gamma_{1k})\cdot g(\gamma_{1k'}).\tag{4.10}$$

Without loss of generality, due to the proof of Lemma 4.4, we can
assume that $\gamma_{1k}=\gamma_{1k_s}\cdot\ \cdots\
\cdot\gamma_{1k_2}\cdot\gamma_{1k_1}$, and that each piecewise
geodesic $\gamma_{1k_j}(t)|_{t\in [0,1]}\subset X_1$ (and
$\tilde\gamma_{1k_j}(t)|_{t\in [0,1]}\subset
(\Sigma_{p_2}X_2)^\perp$) determining  $\gamma_{1k_j}$ satisfies
Case 1 or 2 in that proof. Note that $\gamma_{1k}$ can be determined
(like (4.5)) by the product of $\gamma_{1k_1}(t)|_{t\in
[0,1]},\cdots,\gamma_{1k_s}(t)|_{t\in [0,1]}$ (here, the product
$\alpha\circ\beta(t)|_{t\in[0,1]}$ of curves
$\alpha(t)|_{t\in[0,1]}$ and $\beta(t)|_{t\in[0,1]}$ is defined by
$\alpha\circ\beta(t)=\alpha(2t)$ and $\beta(2t-1)$ for
$t\in[0,\frac12]$ and $[\frac12,1]$ respectively).

Note that $\tilde\gamma_{1k_j}(t)|_{t\in [0,1]}$ is the unique
minimal geodesic between $(\uparrow_{p_2}^{p_1})_1$ and
$(\uparrow_{p_2}^{p_1})_{k_j}$ such that
$f_{p_2}^{-1}(\tilde\gamma_{1k_j}(t)|_{t\in
[0,1]})=\gamma_{1k_j}(t)|_{t\in [0,1]}$ (see (3.2.1)). For
convenience, we say that $\tilde\gamma_{1k_j}(t)|_{t\in [0,1]}$ is
the lifting curve of $\gamma_{1k_j}(t)|_{t\in [0,1]}$ at
$(\uparrow_{p_2}^{p_1})_1$ in $(\Sigma_{p_2}X_2)^\perp$. Note that
$\bar\gamma_{1k_j}((\uparrow_{p_2}^{p_1})_1)=(\uparrow_{p_2}^{p_1})_{k_j},$
so $\bar\gamma_{1k_j}((\uparrow_{p_2}^{p_1})_1)$ is just
$\tilde\gamma_{1k_j}(1)$, the end point of
$\tilde\gamma_{1k_j}(t)|_{t\in [0,1]}$. Furthermore,
$\bar\gamma_{1k_1}(\cdots(\bar\gamma_{1k_s}((\uparrow_{p_2}^{p_1})_1)))$
is the end point of the product of the lifting curves of
$\gamma_{1k_1}(t)|_{t\in [0,1]}$ at $(\uparrow_{p_2}^{p_1})_1$,
$\gamma_{1k_2}(t)|_{t\in [0,1]}$ at
$\bar\gamma_{1k_1}((\uparrow_{p_2}^{p_1})_1)$, $\cdots$, and
$\gamma_{1k_s}(t)|_{t\in [0,1]}$ at
$\bar\gamma_{1k_1}(\cdots(\bar\gamma_{1k_{s-1}}((\uparrow_{p_2}^{p_1})_1)))$
(note that this product is also the lifting curve at
$(\uparrow_{p_2}^{p_1})_1$ of the product of
$\gamma_{1k_1}(t)|_{t\in [0,1]}$, $\cdots$, and
$\gamma_{1k_s}(t)|_{t\in [0,1]}$). This implies that
$$\bar\gamma_{1k_1}\cdot\ \cdots\ \cdot\bar\gamma_{1k_s}=\bar\gamma_{1k}.$$
It follows that
$$g(\gamma_{1k})=\bar\gamma_{1k}^{-1}=\bar\gamma_{1k_s}^{-1}\cdot\
\cdots\ \cdot\bar\gamma_{1k_1}^{-1} =g(\gamma_{1k_s})\cdot\ \cdots\
\cdot g(\gamma_{1k_1}),$$ which implies that (4.10).

\vskip2mm

Due to Claim 1, we let
$\Gamma\triangleq\{\gamma_1,\gamma_2,\cdots,\gamma_m\}$ denote
$\Gamma_1$ and $\Gamma_2$, and let $\Gamma$ act on
$(\Sigma_{p_i}X_i)^\perp$ by isometries such that
$\gamma_k(q_2)=\gamma_{1k}(q_2)$ for any $q_2\in
(\Sigma_{p_1}X_1)^\perp$ and
$\gamma_k(q_1)=\bar\gamma_{1k}^{-1}(q_1)$ for any $q_1\in
(\Sigma_{p_2}X_2)^\perp$. Then the $\Gamma$-action on
$(\Sigma_{p_i}X_i)^\perp$ extends uniquely to an isometric
$\Gamma$-action on
$(\Sigma_{p_1}X_1)^\perp*(\Sigma_{p_2}X_2)^\perp$. Next, through the
following steps, we will construct a 1-1 map
$$\pi:((\Sigma_{p_1}X_1)^\perp*(\Sigma_{p_2}X_2)^\perp)/\Gamma\to X.$$

For convenience, let $\bar X$ denote
$((\Sigma_{p_1}X_1)^\perp*(\Sigma_{p_2}X_2)^\perp)/\Gamma$, and let
$\bar X_j$ denote $(\Sigma_{p_i}X_i)^\perp/\Gamma\subset\bar X$
($j\neq i\in\{1,2\}$), and let $\bar p_j'$ denote any
$\Uparrow_{p_i}^{p_j'}\in\bar X_j$.

\vskip1mm

Step 1:  Define $\pi(\bar p_i')=p_i'$ for any $\bar p_i'\in \bar
X_i\subset\bar X$, i.e., $\pi|_{\bar X_i}=\pi_i.$

Step 2: We define $\pi([\bar p_1\bar p_2'])$ for all $[\bar p_1\bar
p_2']$. Note that there is a natural isometry
$$i_{\bar p_i}:(\Sigma_{\bar p_i}\bar X_i)^\perp\to(\Sigma_{p_i}X_i)^\perp\tag{4.11}$$
(where  $(\Sigma_{\bar p_i}\bar X_i)^\perp$ belongs to $\Sigma_{\bar
p_i}\bar X$). Then we can define $\pi([\bar p_1\bar p_2'])$ to be
$[p_1p_2']$ \footnote{In the present paper, ``$\pi([\bar p_1\bar
p_2'])$ is defined to be $[p_1p_2']$'' means that, for any $\bar
x\in [\bar p_1\bar p_2']$, $\pi(\bar x)$ is defined to be $x\in
[p_1p_2']$ with $|xp_1|=|\bar x\bar p_1|$.} (a minimal geodesics
between $\pi(\bar p_1)$ and $\pi(\bar p_2')$) such that
$\uparrow_{p_1}^{p_2'}=i_{\bar p_1}(\uparrow_{\bar p_1}^{\bar
p_2'})$.

Step 3: For any $\bar p_1'\in \bar X_1$ with $(\pi(\bar p_1')=)\
p_1'\in X_1^m$, we will define $\pi([\bar p_1'\bar p_2'])$ for all
$[\bar p_1'\bar p_2']$ with $\bar p_2'\in \bar X_2$. Select $[\bar
p_1\bar p_1']\subset\bar X_1$ and $[p_1p_1']\subset X_1$ such that
$\pi([\bar p_1\bar p_1'])=[p_1p_1']$. Based on Lemma 4.2 and (4.11),
we can construct an isometry
$$i_{\bar p_1',p_1}\triangleq h_{[p_1p_1']}\circ i_{\bar p_1}\circ h_{[\bar p_1'\bar p_1]}:
(\Sigma_{\bar p_1'}\bar
X_1)^\perp\to(\Sigma_{p_1'}X_1)^\perp\tag{4.12}$$ {\bf Claim 2}:
{\it $i_{\bar p_1',p_1}$ does not depend the choice of
$[p_1p_1']\subset X_1$ (and $[\bar p_1\bar p_1']\subset\bar X_1$)}.
By Remark 4.3, it suffices to verify it for $i_{\bar
p_1',p_1}|_{\{[\bar p_1'\bar p_2]_k\}_{k=1}^m}$, where $\{[\bar
p_1'\bar p_2]_k\}_{k=1}^m$ are all minimal geodesics between $\bar
p_1'$ and $\bar p_2$. In fact, this is clear due to the definition
of $\pi_1\ (=\pi|_{\bar X_1})$ and the isometry $i_{\bar
p_2}:(\Sigma_{\bar p_2}\bar X_2)^\perp\to(\Sigma_{p_2}X_2)^\perp$
(see (4.11)). Based on Claim 2, we can similarly define $\pi([\bar
p_1'\bar p_2'])$ to be $[p_1'p_2']$ such that
$\uparrow_{p_1'}^{p_2'}=i_{\bar p_1',p_1}(\uparrow_{\bar p_1'}^{\bar
p_2'})$.

Step 4: We will define $\pi([\bar p_1'\bar p_2'])$ for all $[\bar
p_1'\bar p_2']$ with $(\pi(\bar p_1')=)\ p_1'\not\in X_1^m$ and
$\bar p_2'\in \bar X_2$. Since $\pi_1$ is an isometry (see (4.9)),
each element of $\Uparrow_{p_2}^{p_1'}\subset
(\Sigma_{p_2}X_2)^\perp$ has an isotropy subgroup $\Gamma_{\bar
p_1'}$ with $|\Gamma_{\bar p_1'}|=\frac{m}{m_{p_1'}}$ (with respect
to the $\Gamma$-action on $(\Sigma_{p_2}X_2)^\perp$). It then
follows that
$$(\Sigma_{\bar p_1'}\bar X_1)^\perp=(\Sigma_{\bar p_1}\bar X_1)^\perp/\Gamma_{\bar p_1'}=
(\Sigma_{p_1}X_1)^\perp/\Gamma_{\bar p_1'}.\tag{4.13}$$ On the other
hand, we have the following cone-neighborhood isometry chain
$$\CD
 X_2 & @>f_{p_1'}>> & (\Sigma_{p_1'}X_1)^\perp & @>f_{\uparrow_{p_1'}^{p_1}}>>
 (\Sigma_{\uparrow_{p_1'}^{p_1}}(\Sigma_{p_1'}X_1))^\perp
\endCD$$
(see (4.1) for $f_{p_1'}$), where $f_{\uparrow_{p_1'}^{p_1}}$ is the
corresponding map to (4.1) for a given $[p_1'p_1]\subset X_1$ in
considering $\Sigma_{p_1'}X\ (\supset \Sigma_{p_1'}X_1,
(\Sigma_{p_1'}X_1)^\perp)$. Then, similar to $h_{[p_1p_1']}$ in
Lemma 4.2, we can define a natural isometry
$$h_{[p_1p_1']}:(\Sigma_{p_1}X_1)^\perp \to
(\Sigma_{\uparrow_{p_1'}^{p_1}}(\Sigma_{p_1'}X_1))^\perp.$$ (Hint:
by Lemma 4.5 (and its proof),
$\{[\uparrow_{p_1'}^{p_1}\uparrow_{p_1'}^{p_2'}]\subset\Sigma_{p_1'}X|\uparrow_{p_1'}^{p_2'}\in
(\Sigma_{p_1'}X_1)^\perp\}$ is mapped 1-1 to $\{[p_1p_2']\subset
X|p_2'\in X_2\}$). Moreover, similar to (4.9), there is a finite
group $\Gamma_{p_1'}$ such that
$$(\Sigma_{p_1'}X_1)^\perp=(\Sigma_{\uparrow_{p_1'}^{p_1}}(\Sigma_{p_1'}X_1))^\perp/\Gamma_{p_1'}
=(\Sigma_{p_1}X_1)^\perp/\Gamma_{p_1'}.\tag{4.14}$$ Based on Lemma
4.6 (and its proof), it is not hard to conclude that $$\Gamma_{\bar
p_1'}=\Gamma_{p_1'}.$$ Therefore, (4.13) and (4.14) enable us to
construct a natural isometry from $(\Sigma_{\bar p_1'}\bar
X_1)^\perp$ to $(\Sigma_{p_1'}X_1)^\perp$, which is similar to
(4.12) and does not depend the choice of $[p_1p_1']$. Then we can
similarly define $\pi([\bar p_1'\bar p_2'])$.

\vskip2mm

So far, we have finished the definition of $\pi$, which is a 1-1
map. We need to show that it is an isometry, i.e.,
$$|\pi(x)\pi(y)|=|xy|
\text{ for all
$x,y\in((\Sigma_{p_1}X_1)^\perp*(\Sigma_{p_2}X_2)^\perp)/\Gamma$}.
\tag{4.15}$$ {\bf Claim 3}: {\it $\pi$ is a continuous map, so it is
a homeomorphism}. It suffices to show that $\pi([\bar p_1^j\bar
p_2^j])\to \pi([\bar p_1'\bar p_2'])$ if $[\bar p_1^j\bar p_2^j]\to
[\bar p_1'\bar p_2']$ as $j\to\infty$. We select $[\bar p_1\bar
p_1^j], [\bar p_1\bar p_1']\subset\bar X_1$ such that $[\bar p_1\bar
p_1^j]\to [\bar p_1\bar p_1']$. By Lemma 1.4, we can find $[\bar
p_1\bar p_2^j]$ and $[\bar p_1\bar p_2']$ such that  $[\bar p_1\bar
p_1^j], [\bar p_1^j\bar p_2^j]$ and $[\bar p_1\bar p_2^j]$ (resp.
$[\bar p_1\bar p_1'], [\bar p_1'\bar p_2']$ and $[\bar p_1\bar
p_2']$) bound a convex spherical surface $D_j$ (resp. $D$), and
$D_j\to D$ (of course, $[\bar p_1\bar p_2^j]\to[\bar p_1\bar
p_2']$). From the definition of $\pi$, it is not hard to see that
$\pi(D_j)$ (resp. $\pi(D)$) is the convex spherical surface bounded
by $\pi([\bar p_1\bar p_1^j]), \pi([\bar p_1^j\bar p_2^j])$ and
$\pi([\bar p_1\bar p_2^j])$ (resp. $\pi([\bar p_1\bar p_1']),
\pi([\bar p_1'\bar p_2'])$ and $\pi([\bar p_1\bar p_2'])$), and that
$\pi([\bar p_1\bar p_1^j])\to\pi([\bar p_1\bar p_1'])$ and
$\pi([\bar p_1\bar p_2^j])\to\pi([\bar p_1\bar p_2'])$. Note that
$\pi(D_j)$ converges to a convex spherical surface including
$\pi([\bar p_1\bar p_1'])$ and  $\pi([\bar p_1\bar p_2'])$. Due to
(4.1.2), it has to hold that $\pi(D_j)\to\pi(D)$ as $j\to \infty$,
so $\pi([\bar p_1^j\bar p_2^j])\to \pi([\bar p_1'\bar p_2'])$.

By Claim 3 and Remark 2.4, for any $[xy]$, $\pi([xy])$ is a
continuous curve with length equal to that of $[xy]$ (note that the
projection of $[xy]$ on $\bar X_i$ is a piecewise geodesic, which is
mapped by $\pi$ to the projection of $\pi([xy])$ on $X_i$), which
implies
$$|\pi(x)\pi(y)|\leq |xy|.$$ Similarly, by considering $\pi^{-1}$, we can conclude that
$$ |xy|\leq|\pi(x)\pi(y)|$$
which implies (4.15).
\hfill$\qed$
\enddemo

%%%%%%%%%%%%%%%%%%%%%%%%%%%%%%%%%%% Section 5%%%%%%%%%%%%%%%%%%%%%%%%%%%%%%%%%%%%%%%%%%%

\head 5. Proof of Corollary 0.2 and 0.3
\endhead

\vskip4mm

\demo{Proof of Corollary 0.2}

Since $M$ and $M_i$ are all Riemannian manifolds, for any $p_i\in
M_i$, each $(\Sigma_{p_i}M_i)^\perp$ is a unit sphere, so is
$(\Sigma_{p_1}M_1)^\perp*(\Sigma_{p_2}M_2)^\perp$. Then Corollary 0.2 is
a corollary of (A2) except when one of $\dim(M_i)$ is equal to 0
(note that in (A2), if $\dim(X_i)=0$, then $X_i$ consists of two points with distance
$\pi$, while in Corollary 0.2, if $\dim(M_i)=0$, then $M_i$ is a single point).
Hence, we need only to consider the case where $\dim(M_1)=0$ (or $\dim(M_2)=0$),
i.e. $M_1=\{p_1\}$.
By Lemma 2.1, we have that $m=1$ or $2$. And since $M$ is a Riemannian manifold,
then it is not hard to see that $M_2^m=M_2$ (see Remark 5.1 below).
By (2.3.2), we can conclude that there is an $m$-order Riemannian cover
from $\Sigma_{p_1}M$ $(=\Bbb S^{n-1})$ to $M_2$. If $m=1$, then $M_2$
is isometric to a unit sphere, so is $M$. If $m=2$, then it is easy
see that $p\in M$ belongs to some $[p_1p_2]$ with $p_2\in M_2$. It then follows
from Lemma 1.4 and (2.3.2) that $\sec_M\equiv1$, so $M$ is isometric to
$\Bbb R\Bbb P^n$.
\hfill$\qed$
\enddemo

\remark{Remark \rm5.1} We would like to point out that in the
case $X$ is a Riemannian manifold, our proof can be made direct
and simple, e.g., discussions in Section 3 and 4 are not required.
Precisely, since $\Sigma_{p_i}M_i$ is a unit sphere, $M_i^m=M_i$,
and a local join structure from (2.3.2) and (2.3.3) implies
that $\sec_M\equiv1$.
\endremark

\vskip2mm

\demo{Proof of Corollary 0.3}

Let $M=\Bbb{S}^n/\Gamma$, where $\Bbb{S}^n$ is the unit sphere. If
$\Gamma=\left<\gamma\right>$, then $\gamma$ is represented by a
closed geodesic on $M/\Gamma$, and $\gamma$ preserves a great circle
$S^1$. Clearly, $\Gamma$ also preserves the sphere
$\Bbb{S}^{n-2}=(S^1)^{=\frac\pi2}$. It then follows that
$M_1=S^1/\Gamma$ and $M_2=\Bbb{S}^{n-2}/\Gamma$ satisfy the
conditions of Corollary 0.2.

In general, $\Gamma$ contains a normal cyclic subgroup $\Bbb Z_q$
with $[\Gamma:\Bbb Z_q]\le w(n)$, a constant depending only on $n$
([Ro]). By the above, the normal covering space of $M$,
$\Bbb{S}^n/\Bbb Z_q$, satisfies Corollary 0.2. \hfill$\qed$
\enddemo

\remark{Remark \rm5.2} From the above proof, we also see the
following: if a spherical space form has diameter less than $\frac
\pi2$, then its fundamental group is not cyclic.
\endremark

%%%%%%%%%%%%%%%%%%%%%%%%%%%%%%%%%%% Section 6%%%%%%%%%%%%%%%%%%%%%%%%%%%%%%%%%%%%%%%%%%%

\head 6. Proof of Theorem B
\endhead

\vskip4mm

The approach to Theorem B is similar to that in the proof of (A2).
The main difference is in showing that for $p_2\in
\partial X_2$,
$f_{p_2}(X_1)=\Uparrow_{p_2}^{X_1}=(\Sigma_{p_2}X_2)^{=\frac \pi2}$
(which guarantees that any $x\in X$ belongs to some $[p_1'p_2']$
with $p_i'\in X_i$), and for any $p_1\in X_1$,
$f_{p_1}(X_2)=\Uparrow_{p_1}^{X_2}= (\Sigma_{p_1}X_1)^{=\frac \pi2}$
(here, $\partial X_2\neq\emptyset$).

For convenience, we list again the main conditions of Theorem B:
$$\dim(X_1)+\dim(X_2)=n-1,\ \text{$X_2^{=\frac \pi2}=X_2^{\ge \frac \pi2}$,\
$(\Sigma_{p_2}X_2)^{=\frac \pi2}
=(\Sigma_{p_2}X_2)^{\ge \frac \pi2}$ for $p_2\in \partial X_2$}.$$

Let's first give a (generalized) version of Key Lemma 1.6 at points
in $\partial X_2$.

\proclaim{Lemma 6.1} For any $p_2\in\partial X_2$, there is an
$m_{p_2}\leq m$ (where $m$ is the number in Key Lemma 1.6) such that
$\lambda_{p_1p_2}\leq m_{p_2}$ for all $p_1\in X_1$, and $\{p_1\in
X_1|\lambda_{p_1p_2}=m_{p_2}\}$ is open and dense in $X_1$.
\endproclaim

\demo{Proof} By Lemma 1.2 and the first variation formula, for a
given $[p_2p_1]$, $|\uparrow_{p_2}^{p_1}\xi|\geq \frac{\pi}2$ for
all $\xi\in \Sigma_{p_2}X_2$. Since
$(\Sigma_{p_2}X_2)^{\geq\frac\pi2}=(\Sigma_{p_2}X_2)^{=\frac\pi2}$,
we in fact have that
$$|\uparrow_{p_2}^{p_1}\xi|=\frac{\pi}2.\tag{6.1}$$
Then for any $[p_2p_2']\subseteq X_2$, there is a $[p_2'p_1]$ which
together with $[p_2p_1]$ and $[p_2p_2']$ bounds a convex spherical
surface (by (1.1.3)). On the other hand, since $X_1$ has empty
boundary, for any $p_2'\in X_2^m$ (ref. Key Lemma 1.6), we still
have that $\Sigma_{p_2'}X=\Sigma_{p_2'}X_2*(\Sigma_{p_2'}X_2)^\perp$
(see (4.1.2)). It then follows that
$$\lambda_{p_1p_2}\leq \lambda_{p_1p_2'}\leq m.$$
Then similar to Lemma 2.6, for any $[p_1p_1']\subset X_1$, we can
conclude that $\lambda_{p_1''p_2}$ is a constant for all $p_1''\in
[p_1p_1']^\circ$, and $\lambda_{p_1p_2}, \lambda_{p_1'p_2}\leq
\lambda_{p_1''p_2}$. This together with Lemma 2.2 implies that there
is an $m_{p_2}\leq m$ such that $\lambda_{p_1p_2}\leq m_{p_2}$ for
all $p_1\in X_1$, and $\{p_1\in X_1|\lambda_{p_1p_2}=m_{p_2}\}$ is
open and dense in $X_1$. (However, so far we do not know that
$X_1^m$ belongs to $\{p_1\in X_1|\lambda_{p_1p_2}=m_{p_2}\}$.)
\hfill$\qed$
\enddemo

Based on Lemma 6.1, using the proof of Key Lemma 1.7, we can
conclude a version of Key Lemma 1.7 at points in $\partial
X_2$.

\proclaim{Lemma 6.2} For any $p_2\in\partial X_2$, $f_{p_2}:X_1\to
(\Sigma_{p_2}X_2)^{=\frac\pi2}$ is also a cone-neighborhood
isometry. As a result, we have that
$\Uparrow_{p_2}^{X_1}=(\Sigma_{p_2}X_2)^{=\frac\pi2}$, which has an
empty boundary.
\endproclaim

\proclaim{Corollary 6.3} For any $x\in X$, there is a $[p_1p_2]$
with $p_i\in X_i$ such that $x\in[p_1p_2]$.
\endproclaim

\demo{Proof} Let $p_2\in X_2$ such that $|xp_2|=\min\{|xp_2'||\
p_2'\in X_2\}$. Note that $|xp_2|\leq\frac\pi2$ because
$X_2^{\geq\frac\pi2}=X_2^{=\frac\pi2}$. By a similar proof of (6.1)
for $p_2\in\partial X_2$ or by Lemma 1.3 for $p_2\in X_2^\circ$, we
know that $\uparrow_{p_2}^x\in (\Sigma_{p_2}X_2)^{=\frac\pi2}$ for
any $[p_2x]$. On the other hand, by Lemma 6.2 or Key Lemma 1.7,
$\Uparrow_{p_2}^{X_1}=(\Sigma_{p_2}X_2)^{=\frac\pi2}$ for any
$p_2\in X_2$. It then follows that there exists $[p_2p_1]$ with
$p_1\in X_1$ such that $\uparrow_{p_2}^x=\uparrow_{p_2}^{p_1}$,
which implies that $[p_2x]\subseteq [p_2p_1]$. \hfill$\qed$
\enddemo

Based on Corollary 6.3, we can give a generalized version of Key
Lemma 1.7 at points in $X_1$ (compared to Theorem A, $\partial
X_2\neq\emptyset$ here).

\proclaim{Lemma 6.4} {\rm (6.4.1)} For any $p_1\in X_1$, $\Uparrow_{p_1}^{X_2}=(\Sigma_{p_1}X_1)^\perp$;
and $f_{p_1}:X_2\to (\Sigma_{p_1}X_1)^\perp$ is also a cone-neighborhood isometry.

\vskip1mm

\noindent{\rm (6.4.2)} For any $p_1\in X_1^m$,
$\Sigma_{p_1}X=\Sigma_{p_1}X_1*(\Sigma_{p_1}X_1)^\perp$.
\endproclaim

\demo{Proof of (6.4.1)}

For the former part, by (1.1) we need only to show that, for any
$\zeta\in(\Sigma_{p_1}X_1)^\perp$, there is a $[p_1p_2]$ with
$p_2\in X_2$ such that $\zeta=\uparrow_{p_1}^{p_2}$. We first assume
that there is a $[p_1x]$ such that $\zeta=\uparrow_{p_1}^x$. Note
that we can assume that there is a unique minimal geodesic between
$x$ and $p_1$ (otherwise we can select $x'\in [p_1x]^\circ$ to take
the place of $x$). By Corollary 6.3, $x$ belongs to some
$[p_1'p_2']$ with $p_i'\in X_i$. It suffices to show that
$p_1'=p_1$. If $p_1'\neq p_1$, then by Lemma 1.4 there is a triangle
$\triangle p_1p_1'p_2'$ containing $[p_1'p_2']$ which bounds a
convex spherical surface. Since there is a unique minimal geodesic
between $x$ and $p_1$, $[p_1x]$ has to lie in this surface. This
implies that $\uparrow_{p_1}^x$ is not perpendicular to the side
$[p_1p_1']$ ($\subset X_1$) in the $\triangle p_1p_1'p_2'$, which
contradicts $\uparrow_{p_1}^x\in (\Sigma_{p_1}X_1)^\perp$. Next, we
need only to show that there must be a $[p_1x]$ such that
$\zeta=\uparrow_{p_1}^x$. If this is not true, we can find
$[p_1x_j]$ such that $\uparrow_{p_1}^{x_j}\to\zeta$ as $j\to\infty$
([BGP]). Using the standard limiting argument we can find a
$[p_1p_2]$ with $p_2\in X_2$ such that $\zeta=\uparrow_{p_1}^{p_2}$,
a contradiction.

For the latter part, we can refer to the proof for (3.9) (here, we
need Lemma 6.1 in addition; and notice that, due to Lemma 6.2, Lemma
1.4 still holds even if $p_2\in \partial X_2$). \hfill$\qed$
\enddemo

Given Lemma 6.2, Corollary 6.3 and (6.4.1), the proof of Lemma 4.5
will go through without a change. That is, we have the following.

\proclaim{Lemma 6.5} In the situation here, the statement in Lemma
4.5 is still true.
\endproclaim

In order to prove (6.4.2), we need the following technical lemma.

\proclaim{Lemma 6.6} Let $A\in \text{Alex}^n(1)$, and let $A_1$,
$A_2$ be two compact convex subsets with $\partial A_1=\emptyset$
and $|A_1A_2|\geq\frac \pi2$. If $\lambda_{a_1a_2}<+\infty$ for all
$a_i\in A_i$ and $\lambda_{a_1a_2}=1$ for $a_2\in A_2^\circ$ and
$a_1\in A_1$, and if any $a\in A$ belongs to some $[a_1a_2]$ with
$a_i\in A_i$, then we have that $A=A_1*A_2$.
\endproclaim

\demo{Proof} By Lemma 1.3, $|a_1a_2|=\frac\pi2$ for all $a_i\in
A_i$, so $A_1*A_2^\circ$ can be isometrically embedded into $A$.
Hence, we need only consider the case where $\partial
A_2\neq\emptyset$ and $\dim(A_2)>0$ (note that if $\dim(A_2)=0$,
then it is our convention that $A_2=A_2^\circ$). {\bf Claim}: {\it
For $a_2\in\partial A_2$, any $[a_2a_1]$ with $a_1\in A_1$ satisfies
$|\uparrow_{a_2}^{a_1}\xi|=\frac\pi2$ for all
$\xi\in\Sigma_{a_2}A_2$}. Then for any $[a_2a_2']\subseteq A_2$ with
$a_2'\in A_2^\circ$, there is a $[a_2'a_1]$ which together with
$[a_2a_1]$ and $[a_2a_2']$ bounds a convex spherical surface (by
(1.1.3)). Since any $a\in A$ belongs to some $[a_1a_2]$ with $a_i\in
A_i$, this implies that
$$A=\overline{A_1*A_2^\circ}, \text{ the closure of } A_1*A_2^\circ$$
(note that $[a_2a_2']\setminus\{a_2\}\subseteq A_2^\circ$). It
therefore follows that $A=A_1*A_2$.

Next, we need only to verify the claim. Consider the multi-valued
map
$$f_{a_2}:A_1\to \Sigma_{a_2}A \text{ defined by } a_1\to \Uparrow_{a_2}^{a_1}.$$
(Note that, by Theorem 1.8, we have that
$\dim(A_1)<n-1=\dim(\Sigma_{a_2}A)$.) Since
$\lambda_{a_1a_2}<+\infty$ for all $a_1\in A_1$, by Lemma 2.2 and
2.6 (like in the proof of Lemma 6.1) we can conclude that there is
an $m$ such that $\lambda_{a_1a_2}\leq m$ for all $a_1\in A_1$ and
$\{a_1\in A_1|\lambda_{a_1a_2}\leq m\}$ is open and dense in $A_1$.
This implies that $f_{a_2}$ is a cone-neighborhood isometry (similar
to (3.9)). Then the claim follows from Lemma 3.8 (note that
$\Sigma_{a_2}A_2\subseteq(f_{a_2}(A_1))^{\geq\frac\pi2}$ in
$\Sigma_{a_2}A$ by the first variation formula). \hfill$\qed$
\enddemo

\demo{Proof of (6.4.2)}

We first consider the case where $\dim(X_2)=0$, i.e., $X_2=\{p_2\}$
because $\partial X_2\ne \emptyset$. By (2.1.1), $m=1$ or $2$.
If $m=1$, then $X$ is isometric to a half suspension (see Corollary 6.3); if $m=2$,
$(\Sigma_{p_1}X_1)^\perp$ consists of two points with distance
$\pi$. In any case, we have that
$\Sigma_{p_1}X=(\Sigma_{p_1}X_1)*(\Sigma_{p_1}X_1)^\perp$.

We now assume that $\dim(X_2)>0$, so
$\dim((\Sigma_{p_1}X_1)^\perp)>0$ by (6.4.1). Note that in
$\Sigma_{p_1}X$, $\partial\Sigma_{p_1}X_1=\emptyset$ and
$\lambda_{\zeta\xi}=1$ for any $\zeta\in
((\Sigma_{p_1}X_1)^\perp)^\circ$ and $\xi\in\Sigma_{p_1}X_1$ (see
(2.11)). In view of Lemma 6.6 with (6.4.1) and Key Lemma 1.7, it
suffices to show $\lambda_{\uparrow_{p_1}^{p_2}\xi}\leq m$ for any
$[p_1p_2]$ with $p_2\in \partial X_2$ and $\xi\in\Sigma_{p_1}X_1$;
and for any $\eta\in \Sigma_{p_1}X$, there is $[\zeta\xi]$ with
$\zeta\in (\Sigma_{p_1}X_1)^\perp$ and $\xi\in\Sigma_{p_1}X_1$ such
that $\eta\in[\zeta\xi]$.

Let $p_2\in \partial X_2$. By Lemma 6.5, 6.1 and 1.4, we can
conclude that for any $[p_1p_1']\subset X_1$, there are at most $m$
minimal geodesics between $\uparrow_{p_1}^{p_2}$ and
$\uparrow_{p_1}^{p_1'}$. I.e., for any $\zeta\in
(\Sigma_{p_1}X_1)^\perp$ (see (6.4.1)), $\lambda_{\zeta\xi}\leq m$
for any $\xi\in(\Sigma_{p_1}X_1)'$, so for any
$\xi\in\Sigma_{p_1}X_1$ by Lemma 2.2.

Next, by Corollary 6.3, any $x\in X$ belongs to some $[p_1''p_2'']$
with $p_i''\in X_i$. Following the argument in the proof of (6.4.1),
we can conclude that $\uparrow_{p_1}^x$ belongs to some
$[\uparrow_{p_1}^{p_1''}\uparrow_{p_1}^{p_2''}]$. This implies that
any $\eta\in \Sigma_{p_1}X$ belongs to some $[\zeta\xi]$ with
$\zeta\in (\Sigma_{p_1}X_1)^\perp$ and $\xi\in\Sigma_{p_1}X_1$.
\hfill$\qed$
\enddemo

\remark{Remark \rm6.7} For $p_2\in \partial X_2$, a priori it is not
clear that $X_1^m\subseteq \{p_1\in X_1|\lambda_{p_1p_2}=m_{p_2}\}$
(see the end of Lemma 6.1). By Lemma 1.4 and (6.4.2), it is not hard
to see that this is true.
\endremark

\vskip2mm

We are now ready to conclude the proof of Theorem B.

\demo{Proof of Theorem B}

In our case where $\partial X_2\ne \emptyset$, the properties from
Lemma 6.1 to Remark 6.7 guarantee that the arguments in the proofs
of Lemma 4.2-4.6 will go through with minor modifications. And the
proof of (A2) also works through here when $\dim(X_2)>0$. So it
remains to consider the case where $\dim(X_2)=0$, i.e.,
$X_2=\{p_2\}$.

By (2.1.1), we have that $\lambda_{p_1p_2}\leq m$ for all $p_1\in
X_1$  with $m=1$ or 2. If $m=1$, then $X$ is a half suspension of
$X_1$ (see Corollary 6.3). If $m=2$, then there is an obvious
involution on $\Sigma_{p_2}X$, which is an isometry (similar to
$\hat \gamma$ in (4.8)), and on
$(\Sigma_{p_1}X_1)^\perp=\{(\uparrow_{p_1}^{p_2})_1,
(\uparrow_{p_1}^{p_2})_2\}$ for any $p_1\in X_1^m$. Then it is not
hard to see that $X=(\Sigma_{p_2}X*(\Sigma_{p_1}X_1)^\perp)/\Bbb
Z_2=S(\Sigma_{p_2}X)/\Bbb Z_2$ (and $X_1=(\Sigma_{p_2}X)/\Bbb Z_2$).\hfill$\qed$
\enddemo

\enddocument